\documentclass[12pt]{article}

\newcommand{\subfigimg}[3][,]{%
	\setbox1=\hbox{\includegraphics[#1]{#3}}
	\leavevmode\rlap{\usebox1}
	\rlap{\hspace*{-5pt}\raisebox{\dimexpr\ht1-2\baselineskip}{#2}}
	\phantom{\usebox1}
}

\usepackage[utf8]{inputenc}

\usepackage{authblk}
\usepackage{amsmath}
\usepackage{amsfonts}
\usepackage{amssymb}
\usepackage{graphicx}
\usepackage{color}
\usepackage{hyperref}
\usepackage[margin=1in]{geometry}

\graphicspath{{./img/}}

\newtheorem{theorem}{Theorem}[section]
\newtheorem{remark}[theorem]{Remark}

\usepackage{soul}

\newcommand{\RED}[1]{{\color{black}{#1}}}
\newcommand{\yw}[1]{{\color{magenta}{#1}}}

\title{A conceptual framework for modeling a latching mechanism for cell cycle regulation}
\author[1]{Punit Gandhi} 
\author[2]{Yangyang Wang}
\affil[1]{\small{Department of Mathematics and Applied Mathematics, Virginia Commonwealth University, Richmond, Virginia, USA 23284}}
\affil[2]{\small{Department of Mathematics, Brandeis University, Waltham, Massachusetts, USA, 02453}}
\date{}

\begin{document}
\maketitle

\begin{abstract}
Two identical van der Pol oscillators with mutual inhibition are considered as a conceptual framework for modeling a latching mechanism for cell cycle regulation.  In particular, the oscillators are biased to a latched state in which there is a globally attracting steady-state equilibrium without coupling.  The inhibitory coupling induces stable alternating large-amplitude oscillations that model the normal cell cycle.  A homoclinic bifurcation within the model is found to be responsible for the transition from normal cell cycling to endocycles in which only one of the two oscillators undergoes large-amplitude oscillations.  
\end{abstract}





\section{Introduction}

\begin{figure}[htp!]
    \begin{center}
\begin{tabular}
{@{}p{0.33\linewidth}@{\quad}p{0.33\linewidth}@{\quad}p{0.33\linewidth}@{}}
\subfigimg[width=\linewidth]{\bfseries{\small{(A)}}}{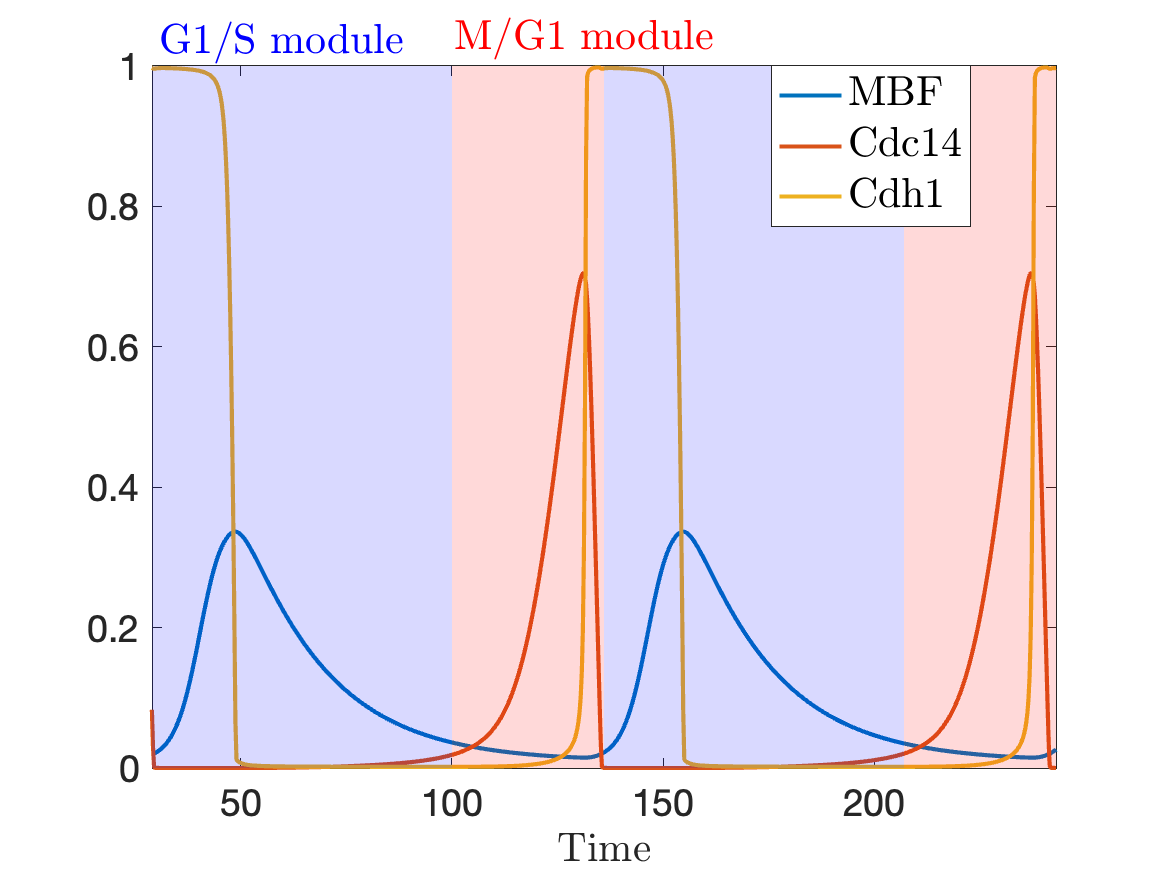}& 
   \subfigimg[width=\linewidth]{\bfseries{\small{(B)}}}{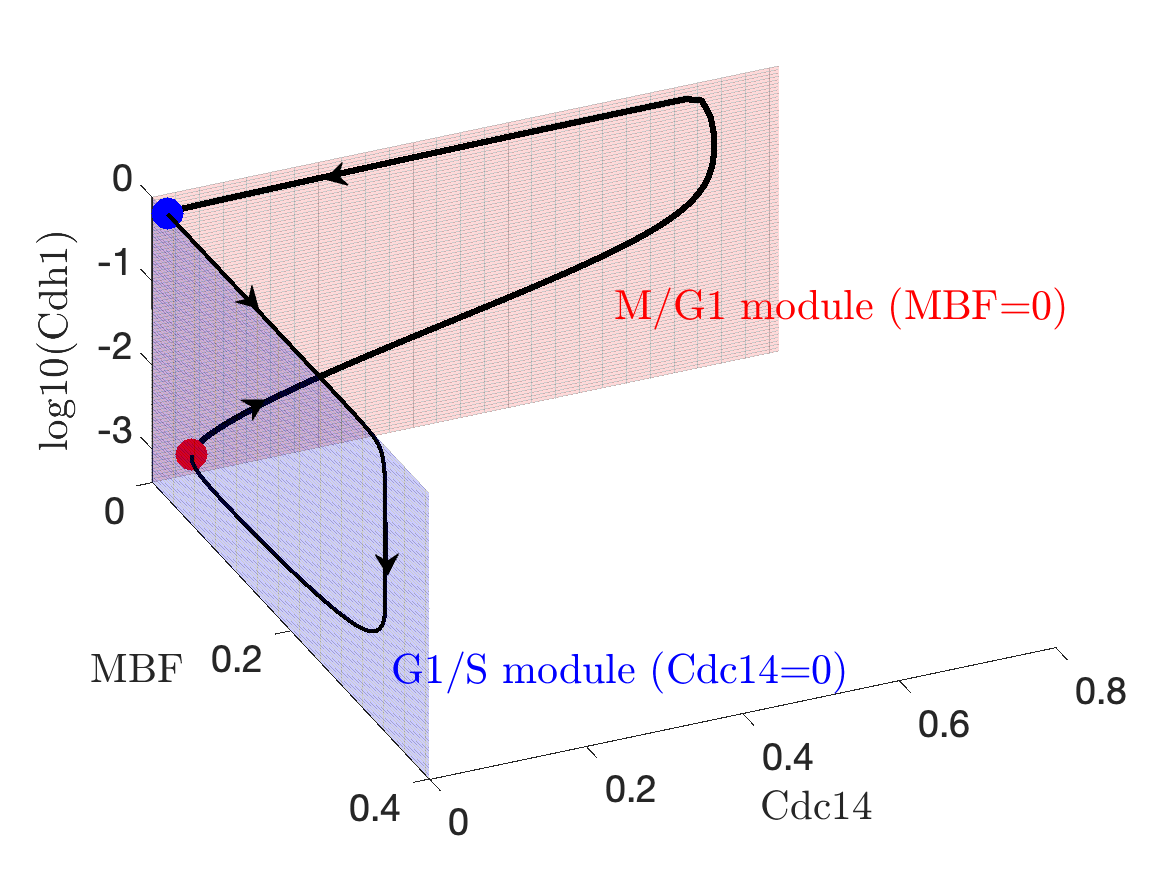}&
    \subfigimg[width=\linewidth]{\bfseries{\small{(C)}}}{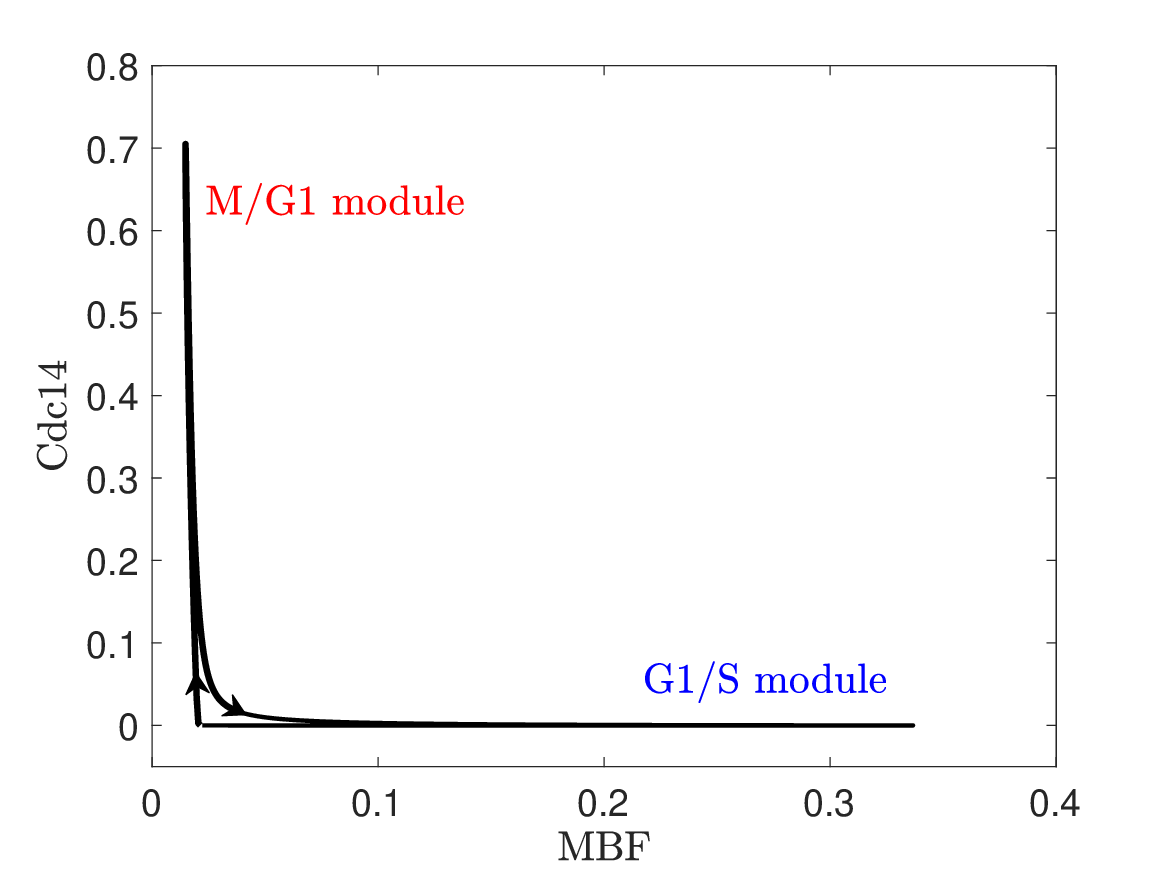}
\end{tabular}
\end{center}
 \caption{The normal cell division cycle in wild-type budding yeast. The full model and parameters are given in \cite{novak2022mitotic}. (A) Simulated time course of a mother cell undergoing two divisions. (B,C) Projections of the solutions from (A) onto $(\rm MBF,Cdc14,log10(Cdh1))$ and $(\rm MBF,Cdc14)$, respectively.}
 \label{fig:biomodel}
\end{figure}

The eukaryotic cell cycle is a recurrent alternating sequence of ``one-way" transitions from \RED{G1} (unreplicated chromosomes) to \RED{S-G2-M} (replication and partitioning of chromosomes), and back to \RED{G1} again \cite{nasmyth1996,tyson1995,novak2007,tyson2008,novak2022mitotic,Dragoi2024cellcycle,Dragoi2024model}. \RED{Novak} and Tyson~\cite{novak2022mitotic} propose a latching-gate mechanism for regulation in the yeast cell cycle, where the normal mitotic cycles can be viewed as a bistable switch between the two stable steady states \RED{G1} (APC/C:Cdh1 activity high) and \RED{S-G2-M} (APC/C:Cdh1 activity low). These two phases correlate to high activity of cyclin-dependent kinases (CDK) that trigger \RED{S-G2-M} phase and CDK antagonists (e.g., Cdh1) that stabilize the $G_1$ phase. Figure \ref{fig:biomodel} shows the numerical simulation results of the cell cycle latching model described in~\cite{novak2022mitotic}, depicting the normal cell division cycle in wild-type budding yeast. 

In order to enter \RED{S-G2-M} from the \RED{G1} phase (Figure \ref{fig:biomodel}A, \RED{G1/S} module in blue shaded region), the gate must be sufficiently pushed in one direction by a transient activation of the CDK transcription factor `MBF' (denoted by the blue curve). Once this activation occurs, the gate swings back into its latched position. To return to the \RED{G1} phase (Figure \ref{fig:biomodel}A, \RED{M/G1} module in the red shaded region), the gate must be pushed in the opposite direction by a transient activation of a mitotic exit protein that activates Cdh1 (`Cdc14', represented by the red curve), after which the gate latches again.

Specifically, starting at a pseudo steady state (the \RED{G1} steady state denoted by the blue circle in Figure \ref{fig:biomodel}B), the transient activation of MBF is needed to open the gate, allowing the cell to leave the red plane and progress from \RED{G1} to \RED{S-G2-M} (illustrated by the trajectory along the blue plane). During this transition, negative feedback on MBF pulls the gate closed again, and the ``latch" (a pseudo steady state at the red circle) holds the cell in mitosis. To exit mitosis and return to the \RED{G1} phase, Cdc14 must be activated (see the trajectory near the red plane). After the activation, the negative feedback on Cdc14 brings the trajectory back to the blue circle, and the latch holds the cell in \RED{G1} until MBF is reactivated, completing a full cycle. The projection of the solution trajectory onto $(\rm MBF, Cdc14)$-space is shown in Figure \ref{fig:biomodel}C, revealing an interesting butterfly shape with two wings or loops that represent the alternation between the two modules. Moreover, MBF and Cdc14 act as ``latching parameters" that unlatch the gates so that Cdh1 can cycle from low to high and back again. \cite{novak2022mitotic} shows that perturbations that break the latch may convert a latching gate mechanism into a swinging gate and produce either \RED{G1-S-G1-S} (periodic DNA replication without mitosis \cite{edgar2001}, denoted as ``endoreplication") or \RED{M-(G1)-M-(G1)} (periodic Cdc14 release without fully exiting mitosis \cite{lu2010,manzoni2010}, denoted as ``Cdc14 endocycles"). Such one-loop oscillations are called \emph{endocycles}.

The latching model for the cell cycle progression in \cite{novak2022mitotic} consists of a system of 10 differential equations, making it difficult to identify the key regulatory mechanism driving the latching dynamics. Although bifurcation analysis has been conducted, it was performed on reduced systems, assuming some dynamic variables were zero, rather than on the full system.
In this paper, we propose a simplified conceptual model 
that captures both normal mitotic cycles, characterized by a strict alternation of \RED{G1/S} and \RED{M/G1} modules as shown in Figure \ref{fig:biomodel}, and endocycles, which span only one of these two modules. 
We require (1) a transient activation of Cdc14 to unlatch the gate for Cdh1 to go from low to high and (2) a pulse of MBF to unlatch the gate for Cdh1 to go from high to low (see Figure \ref{fig:biomodel}B). To achieve this, we propose using coupled van der Pol oscillators \cite{vanderpol1920,vanderpol1926}: one to regulate Cdc14 oscillations and the other to control MBF oscillations. 



Inspired by efforts to develop a simple conceptual model that captures the mathematical mechanism underlying the cell cycle control in budding yeast \cite{novak2022mitotic}, we consider a pair of identical coupled van der Pol oscillators of the form
\begin{subequations}\label{eq:cVdP}
\begin{align}
     \dot{x}_1 &=y_1 + \left(1 -\frac{x_1^2}{3}\right)x_1\\
    \dot{y}_1 &=\varepsilon_1(a_1-x_1 + b_{1} \tanh(k_2(a_2-x_2))) \\ 
    \dot{x}_2 &=y_2 + \left(1 -\frac{x_2^2}{3}\right)x_2\\
    \dot{y}_2 &=\varepsilon_2 (a_2-x_2 + b_{2}\tanh(k_1(a_1-x_1))) 
 \end{align}
\end{subequations}
with timescale separations $\epsilon_i$, latching parameter $a_i$ and coupling strengths $b_{i}$. The dot represents the derivative with respect to time $t$. Throughout the paper, we consider identical oscillators by assuming $\varepsilon_1=\varepsilon_2=\varepsilon$, $a_1=a_2=a$, $b_1=b_2=b$ and $k_1=k_2=k$. All the parameters are given in Table \ref{tab:parameters}. Without loss of generality, we assume $a\leq 0$, as the dynamics for $a\geq 0$ will be symmetric to dynamics for $a\leq 0$. 

In this model, $(x_1,y_1)$ represents the \RED{G1/S} module with $x_1$ corresponding to MBF, while $(x_2,y_2)$ represents the \RED{M/G1} module with $x_2$ corresponding to Cdc14.  \RED{We emphasize that $(x_i,y_i)$ are intended to capture qualitative trends of the corresponding biological variables and their values should not be taken as concentrations. In particular, they can take on both positive and negative values.  We choose the form of the nonlinear coupling terms in~\eqref{eq:cVdP} in order to qualitatively match the  latched state dynamics described in~\cite{novak2022mitotic}. A positive coupling strength $b$ implies mutual inhibition between the two modules, and the $\tanh$ functional form  ensures saturation of these inhibitory effects.  The difference $y_1-y_2$ serves as an analog for Cdh1, as it reproduces the qualitative behavior shown in Figure~\ref{fig:biomodel}. This analogy aligns qualitatively with the dynamics of the cell cycle model in~\cite{novak2022mitotic}, where $x_1$ inhibits $y_1$ (analogous to MBF inhibiting Cdh1), and $x_2$ activates $-y_2$ (analogous to Cdc14 activating Cdh1).  While $y_1$ and $-y_2$ individually capture aspects of the inhibition and activation dynamics of Cdh1, it is the combination $y_1-y_2$ that reflects the overall behavior of Cdh1 in the model. 
} 

For any coupling strength \RED{$b$}, system \eqref{eq:cVdP} has a \emph{symmetric steady-state solution} (denoted as $E_0$) that satisfies 
\begin{equation}\label{eq:eq}
    x_1=x_2=a, \; y_1=y_2=-(1-a^2/3)a.
\end{equation}
With parameters in Table \ref{tab:parameters}, $E_0$ is given by 
$x_1 = x_2 = -1.1, y_1 = y_2 = 0.6563$. 
To capture the butterfly-like dynamics from the biological model (see Figure \ref{fig:biomodel}C), where during each module either $x_1$ (MBF) or $x_2$ (Cdc14) is maintained at a constant level near its pseudo-steady state, we selected parameters such that the uncoupled oscillators exhibit stable equilibrium (i.e., the latched state) when isolated. We show later in the Results section \ref{sec:result} that the two oscillators can interact through mutual inhibition to switch back and forth between the two modules, as observed in normal cell cycle progression. We also show that increasing the inhibitory coupling strength $b$ can break the latch mechanism in \eqref{eq:cVdP} and induce the transition from the normal double-loop oscillations to single-loop endocycles. 


\begin{table}[!htb]
	\renewcommand{\arraystretch}{1.55}
	\centering
    \begin{tabular}{|c|c|c|l|}
    \hline
         & $i=1,2$  & description\\
         \hline
        $\epsilon $ & $10^{-2}$  & timescale separation  \\
         $a$ & -1.1   & latching parameter \\
         $b$ & 0.3   & inhibitory coupling strength \\
         $k$ & 1 & coupling response sensitivity\\
         \hline
    \end{tabular}
    \caption{Parameters for a pair of identical coupled van der Pol oscillators in Eq.~\eqref{eq:cVdP}.}
    \label{tab:parameters}
\end{table}


Analysis of our conceptual model suggests that a homoclinic bifurcation  plays a key role in the transition to endocycles within our model~\eqref{eq:cVdP} of two identical coupled oscillators.  We see a pair of periodic orbits collide with a steady-state equilibrium in what has been called a \textit{homoclinic butterfly} in the Lorenz equations (See, e.g., \cite{shilnikovL2003HC}).  Its role in the path to chaotic dynamics has been explored within the Lorenz equations and other related models~\cite{sparrow1982lorenz,xing2014symbolic,pusuluri2021homoclinic}. Homoclinc bifurcations with $\mathbb{Z}_2$ symmetry have also been studied in the context of models of Josephson junction arrays~\cite{aronson1991HC,aronson1994HC}.  
 
\RED{Our work represents a conceptual framework for capturing swinging and latching dynamics in a mathematically tractable system.  The structure underlying model~\eqref{eq:cVdP} is rooted in the idea of \textit{echo waves}~\cite{krinskii1972investigation,winfree1972spiral,tyson1979oscillations} wherein coupled excitable media generate large-amplitude, out-of-phase oscillations. 
Model~\eqref{eq:cVdP} does not attempt to capture details of the biological processes involved, but instead takes a phenomenological approach for qualitatively reproducing the observed behavior in~\cite{novak2022mitotic}. 
The related work of Dragoi et al. \cite{Dragoi2024cellcycle} also uses identical coupled oscillators to model the latching dynamics in a complex cell cycle model, and remains closely connected to the known biochemistry of cell cycle regulation in eukaryotes.  A key biological phenomenon of interest in our study is the transition from normal cell cycling to endocycles, which is mathematically represented by a transition from double-loop to single-loop oscillations.  Interestingly, the mechanism responsible for this transition in the model~\eqref{eq:cVdP} is different from what is observed by Dragoi et al.~\cite{Dragoi2024cellcycle}, which we discuss further in Section~\ref{sec:discuss}.}  

The transition from double-loop to single-loop oscillations, recently referred to as \emph{strong symmetry breaking}, have \RED{also} been examined by \cite{awal2023symmetry,awal2024symmetry} in coupled identical oscillators. These studies highlight the critical role of specific folded singularities near the symmetric equilibrium driving the strong symmetry breaking. In contrast, \RED{model~\eqref{eq:cVdP}} 
does not rely on folded singularities to shape its dynamics; in fact, folded nodes are not present within the considered parameter regime (see \ref{app:gspt} for more details). A detailed comparison between our work and other related studies will be presented in Section \ref{sec:discuss}. 

Our paper is organized as follows. The main results are presented in Section \ref{sec:result}, which is divided into four subsections. 
We begin by analyzing the dynamics of single uncoupled oscillator in Subsection \ref{sec:single} and demonstrate that $x_i$ can be interpreted as ``latching" parameters, consistent with the biological model. In Subsection \ref{sec:simu-coup}, we show that while the uncoupled oscillators remain latched in isolation, a sufficient push of the system through mutual inhibition can momentarily open the latched gate, producing oscillations that alternate between the \RED{G1/S} and \RED{M/G1} modules (i.e., double-loop oscillations). By symmetrically increasing the inhibitory coupling strength $b$ in \eqref{eq:cVdP}, we observe a strong symmetry-breaking phenomena similar to that described in \cite{awal2024symmetry,epstein2024}, leading to the transition from double-loop oscillations to single-loop endocycles. To gain deeper insight into the mechanism driving this transition, we perform a one parameter bifurcation analysis of the coupled-oscillator system~\eqref{eq:cVdP} by treating $b$ as the bifurcation parameter in Subsection \ref{sec:onepara}. This analysis reveals a homoclinic (HC) bifurcation that separates the double-loop and single-loop solutions. We demonstrate in Subsection \ref{sec:HCButterfly} that this HC bifurcation is the key mechanism driving the strong symmetry breaking in our model. In Section \eqref{sec:twopara}, we further investigate the key bifurcations by varying an additional latching parameter, $a$. We demonstrate that the HC bifurcation aligns well with the boundary between the double-loop and single-loop oscillations in $(a,b)$ parameter space. 
Finally, Section \ref{sec:discuss} contains a discussion of our results and future work.  

\section{Results}\label{sec:result}
\subsection{Single oscillator $(x_i, y_i)$}\label{sec:single}

\begin{figure}[!htp]
    \centering   
    \includegraphics[width=3in]{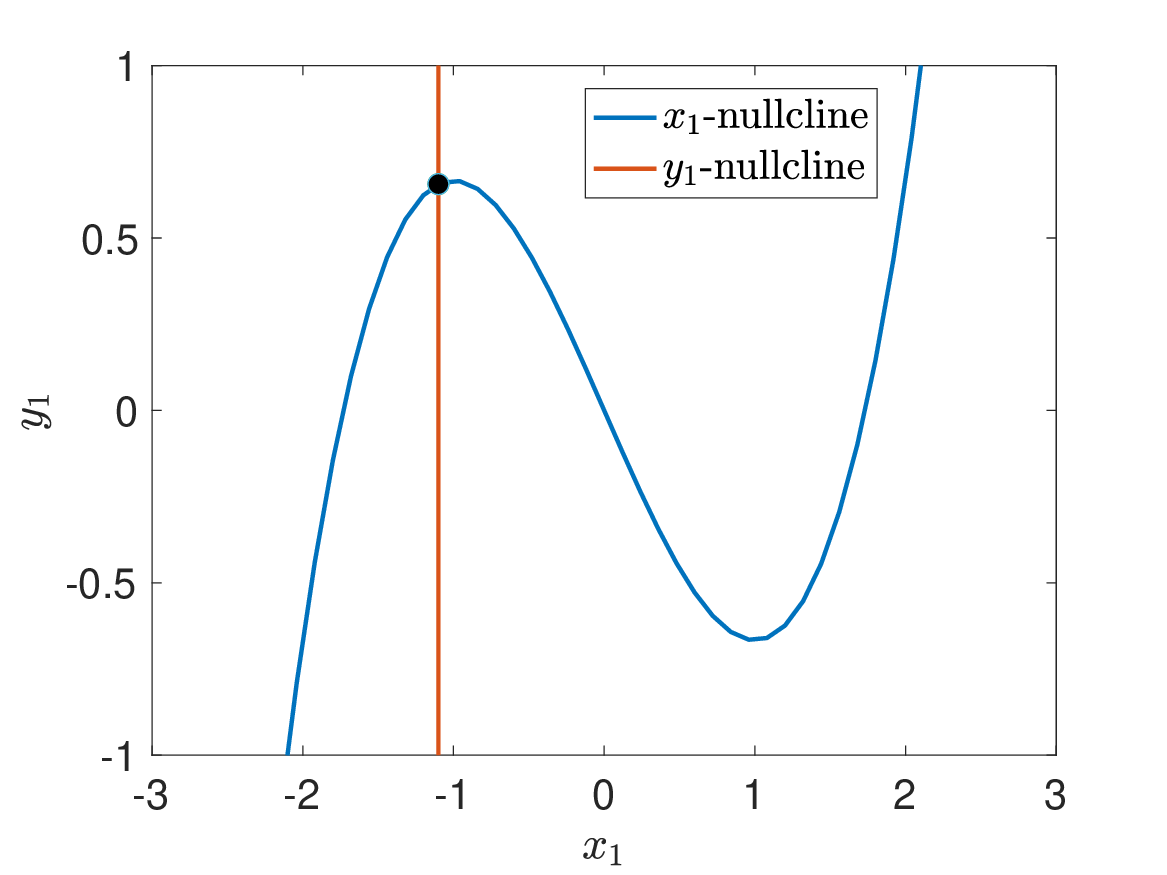}
    \includegraphics[width=3in]{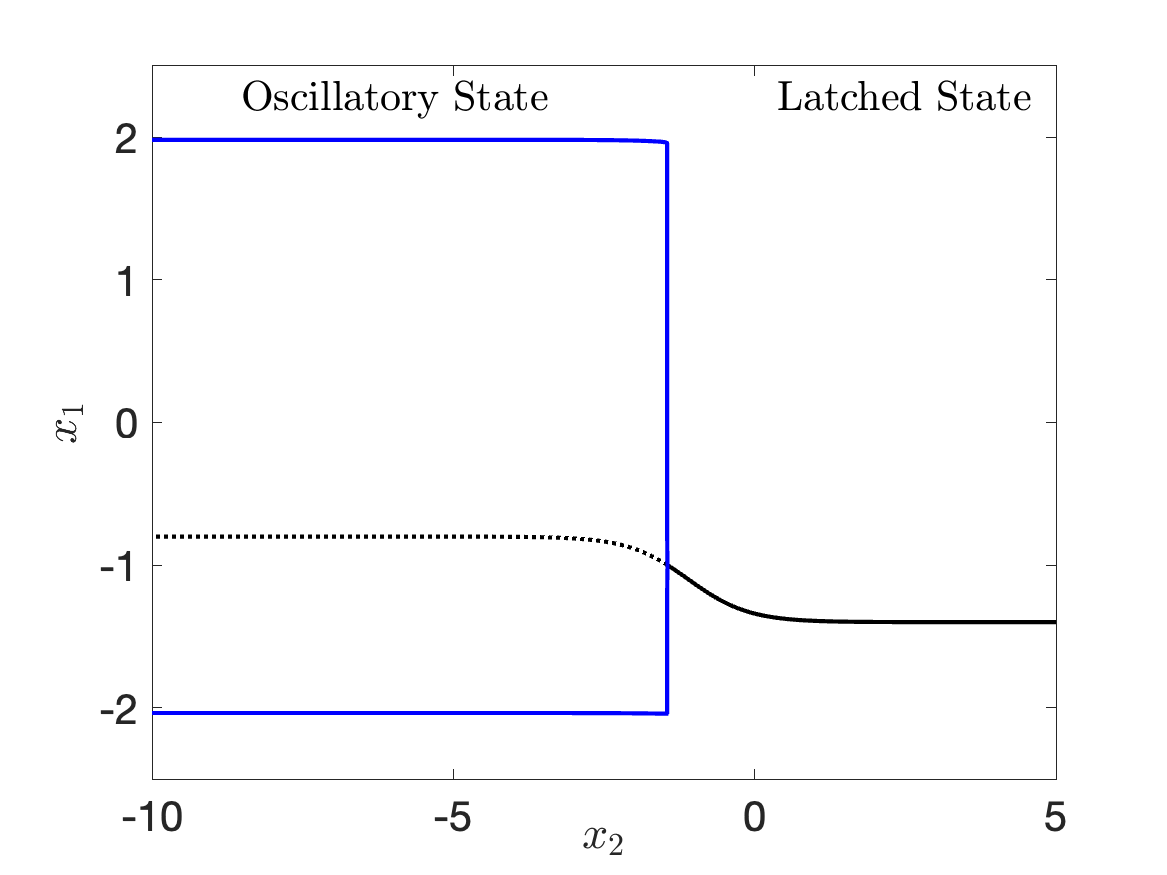}
    \caption{Dynamics of the single oscillator with parameters are given in Table \ref{tab:parameters}. Left:
    Phase portrait of the single oscillator $(x_1, y_1$) without coupling $b=0$. The two nullclines intersect at a stable fixed point at the black circle. Right: With coupling $b=0.3$, $x_2$ serves as the ``latching" parameter for the oscillator $(x_1, y_1$): When $x_2$ is high, the oscillator is being put in the latched state ($y_1$-nullcline is intersecting $x_1$-nullcline at its left branch); As $x_2$ decreases, the oscillator $(x_1,y_1)$ moves to the oscillatory state ($y_1$-nullcline moves rightward to intersect $x_1$-nullcline at its middle branch). The superscript HB bifurcaiton occurs at $x_2 = -1.447$. }
    \label{fig:single_vdp}
\end{figure}

When uncoupled ($b=0$), the Jacobian matrix of the single oscillator $(x_i,y_i)$ at its equilibrium $(a_i, (-1+\frac{a_i^2}{3})a_i)$ is given by 
\[
J_i = \begin{bmatrix}
  (1-a_i^2) & 1\\
    -\varepsilon_i & 0
\end{bmatrix}. 
\]
The Hopf bifurcation points are given by ${\rm tr}(J)= (1-a_i^2) = 0$ and $\det J=\varepsilon_i>0$. It follows that the single oscillator exhibits a Hopf bifurcation at \RED{$a_{HB}=- 1$} (we omit the other HB value at $a_i=1$ as discussed earlier). In \eqref{eq:cVdP}, we choose \RED{$a_i=-1.1<a_{HB}$} so that the oscillators are biased just beyond the Hopf bifurcation for the uncoupled oscillators so that the two nullclines intersect at a globally attracting fixed point (see Figure \ref{fig:single_vdp}, left). \RED{We call this stable steady-state the latched state, because the system has no limit cycle but perturbations can induce a single large-amplitude oscillation cycle before returning to this latched state. Increasing the parameter $a_i$ across $a_{HB}=-1$ is able to switch the system from the latched state to the swinging state in which a stable limit cycle solution exists, leading to sustained periodic oscillations. Thus, we refer to $a_i$ as a \emph{latching parameter}.}

The right panel of Figure \ref{fig:single_vdp} shows the effect of the inhibition of $x_2$ \RED{on the first oscillator $(x_1,y_1)$}. Decreasing $x_2$, which shifts the $y_1$-nullcline to the right, can transition the system from a latched state to an oscillatory state after passing a supercritical HB bifurcation at $x_2=-1.447$. For $x_2<-1.447$, the two nullclines of $(x_1,y_1)$ intersect at an unstable steady-state in the middle branch of the $x_1$-nullcline, and the oscillator exhibits limit cycle oscillations. 

\RED{For our choice of $a_i<a_{HB}$,} these uncoupled oscillators are excitable in the sense that sufficiently large perturbations away from the fixed point can lead to large excursion of the system before eventually returning to the fixed point. In other words, a sufficient push of the system can momentarily open the latched gate, but it will not continue oscillating when uncoupled. This behavior differs from other reduced cell cycle model\RED{s} (e.g., \cite{Dragoi2024cellcycle}), where each oscillator exhibits stable limit cycle oscillations even in isolation.

\subsection{Double-loop oscillations with inhibitory coupling}\label{sec:simu-coup}

\begin{figure}[!htp]
\begin{center}
\begin{tabular}
{@{}p{0.5\linewidth}@{\quad}p{0.5\linewidth}@{}}
\subfigimg[width=\linewidth]{\bfseries{\small{(A)}}}{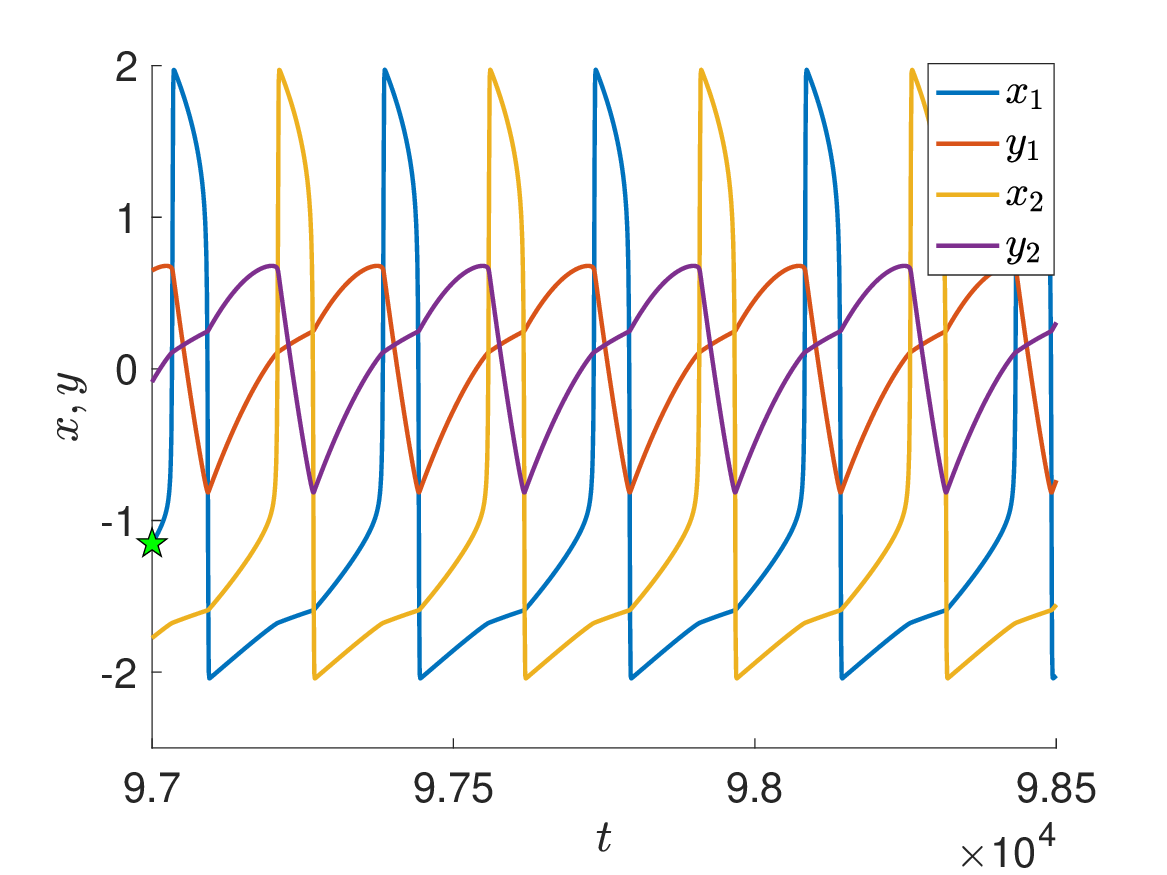}&
\subfigimg[width=\linewidth]{\bfseries{\small{(B)}}}{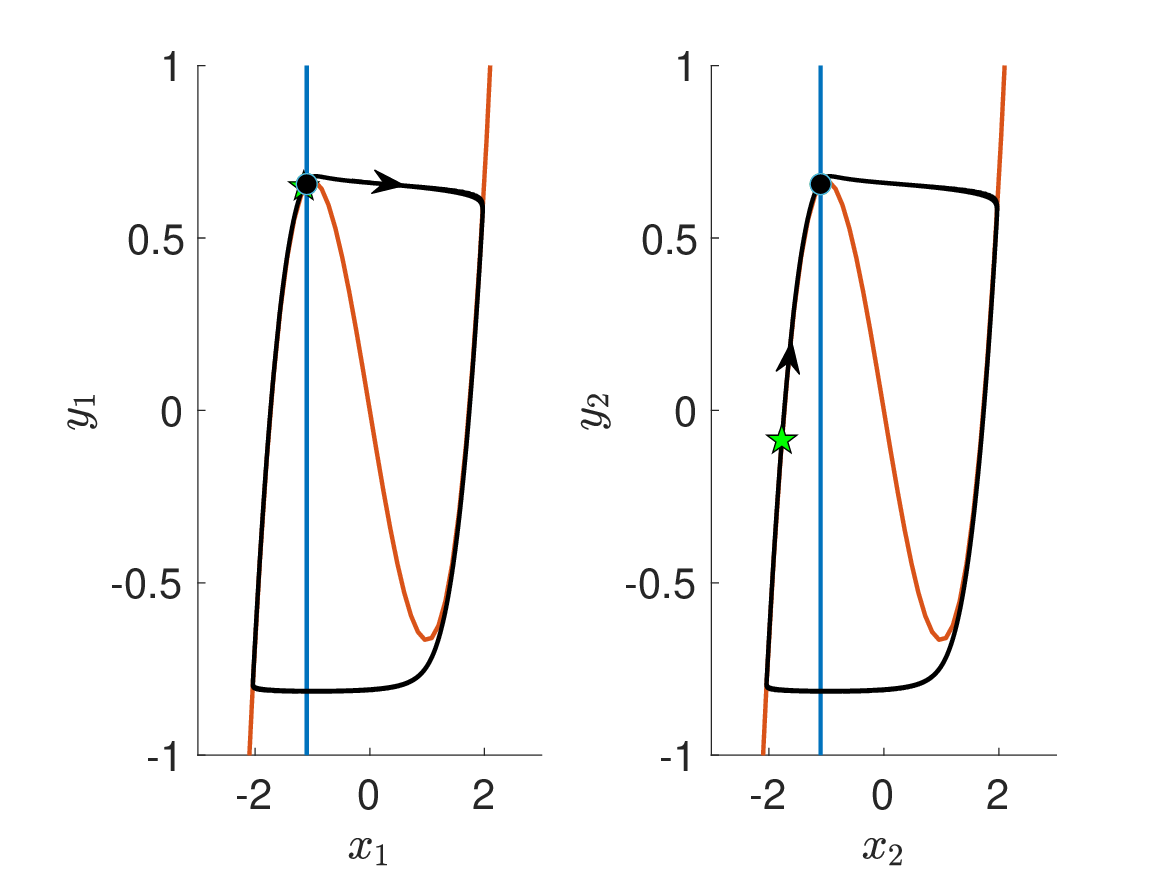}\\
\subfigimg[width=\linewidth]{\bfseries{\small{(C)}}}{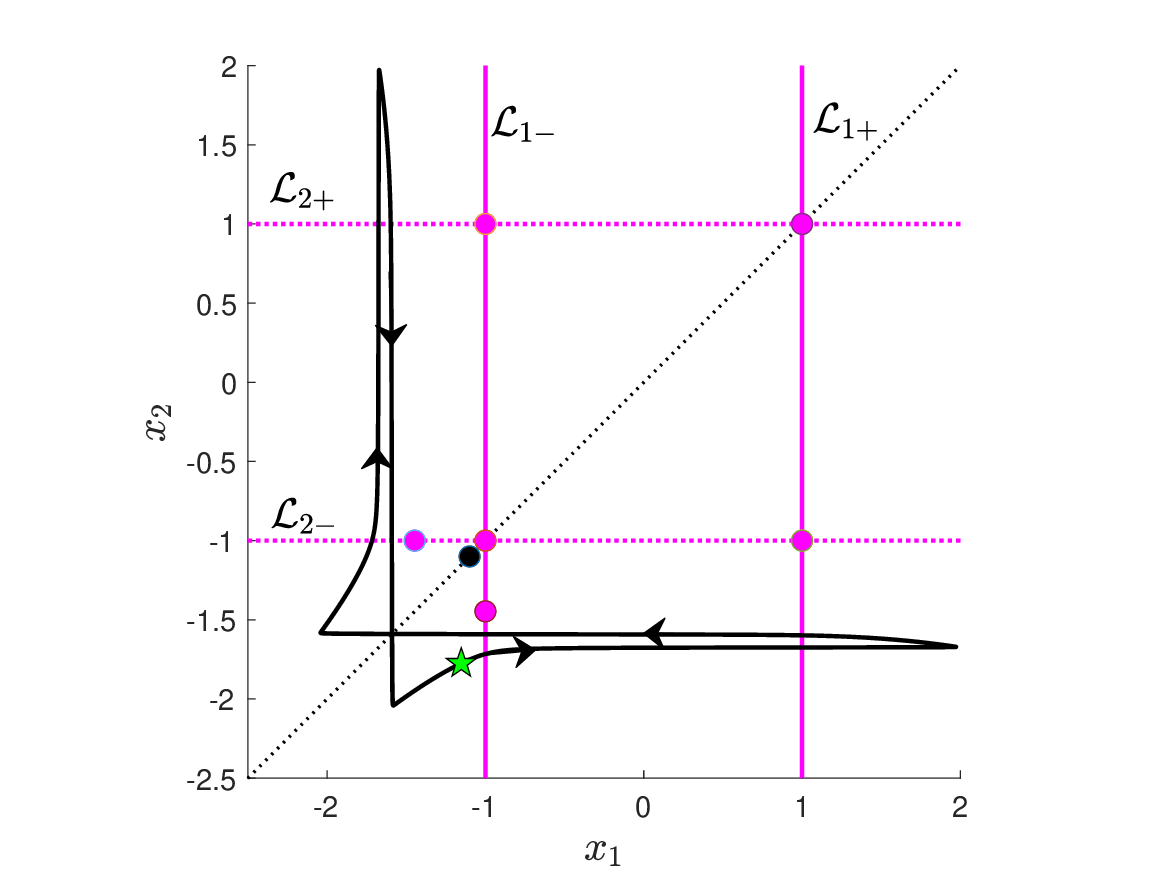} & \subfigimg[width=\linewidth]{\bfseries{\small{(D)}}}{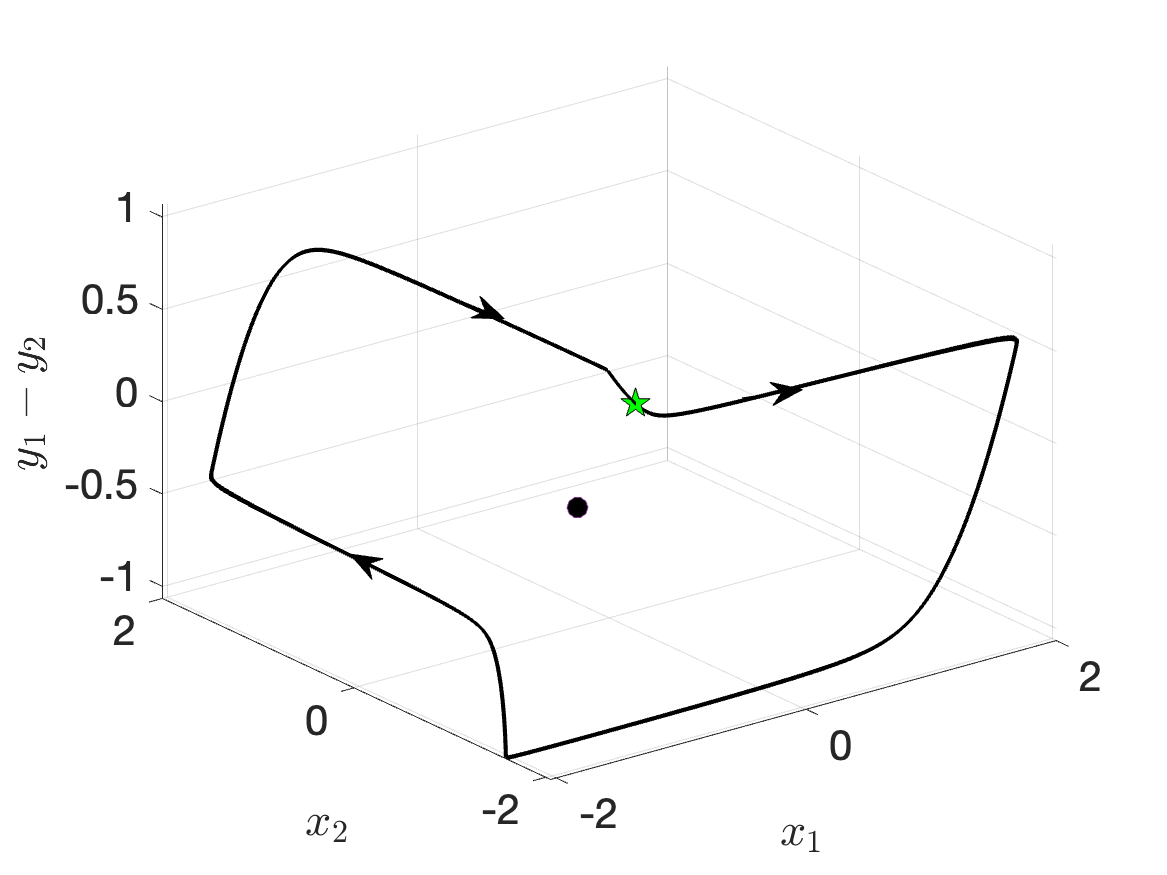}
\end{tabular}
\end{center}
\caption{Solutions of \eqref{eq:cVdP} at the oscillatory state for default parameters given in Table \ref{tab:parameters}. (A) Time traces of $(x_i, y_i)$. (B-D) Projections of solutions from (A) onto $(x_i,y_i)$-space, $(x_1,x_2)$-space and $(x_1,x_2,y_1-y_2)$-space. Green star marks the starting position. Black circle indicates the true equilibrium ($x_i=1.1$, $y_i=0.6563$), which is independent of the coupling strength. There is a bistability between this stable equilibrium and the oscillation (black lines). Red and blue lines in panel (B) indicate $x_i$ and $y_i$-nullclines, respectively, when uncoupled. In panel (C), magenta lines indicate folds of the $x_i$-nullclines, while the magenta circles on these fold curves denote folded singularities. See Figure \ref{fig:fs} for the folded singularity types.} 
\label{fig:soln}
\end{figure}

Figure \ref{fig:soln}A shows the time-course simulations of \eqref{eq:cVdP} with inhibitory coupling strength $b=0.3$. The two systems oscillate out of phase, similar to the out-of-phase oscillations observed between the \RED{G1/S}- and \RED{M/G1}- modules in the cell cycle model \cite{novak2022mitotic}.
Although each system is in the latched state (denoted by the black circle in Figure \ref{fig:soln}) in isolation, when the coupling is introduced, oscillations are possible because one latched gate can be opened during the last phase of the large excursion of the other (See Figure \ref{fig:soln}B). As previously discussed, the black circle is the symmetric equilibrium $E_0$ for the coupled system for any coupling strength (see \eqref{eq:eq}). When $b=0.3$, system \eqref{eq:cVdP} exhibits bistability between $E_0$ and a stable limit cycle.

Next, we explain how mutual inhibition between the two oscillators generates double-loop oscillations, even when each system independently has a stable latched state. Starting at the green star, $x_2$ is sufficiently small (below the HB bifurcation point in Figure \ref{fig:single_vdp}, right), putting the first oscillator in an oscillatory state. Consequently, the trajectory in $(x_1,y_1)$ jumps to the right branch of $x_1$-nullcline to make a large excursion (c.f., \RED{G1/S} module in Figure \ref{fig:biomodel}). After passing the lower fold of the $x_1$-nullcline, the trajectory jumps back to the left branch, where $x_1$ now falls below the HB bifurcation, unlatching the second oscillator and allowing it to swing freely (c.f., \RED{M/G1} module). This process continues until the trajectory returns to the green star, completing a full oscillatory cycle. In each cycle, the transient activation of $x_1$ or $x_2$ that unlatches the gate occurs as $x_2$ or $x_1$ falls below the Hopf bifurcation.  

Figure \ref{fig:soln}C and D show the projections of the solutions onto $(x_1,x_2)$-space and $(x_1,x_2,y_1-y_2)$-space, respectively, analogous to the representations of the biological model dynamics shown in Figure \ref{fig:biomodel}C and B. The similarities between these projections further validate the ability of our model to capture the key dynamics of the cell cycle model from \cite{novak2022mitotic}. 
In panel (C), we also plot the four fold lines of the $x_i$-nullclines given by \eqref{eq:fold} and a total of six folded singularities described in \ref{rm:fs}. Unlike the findings in \cite{awal2023symmetry}, we observe only folded saddles and folded foci near the symmetric equilibrium $E_0$, which do not play a significant role in shaping the dynamics in \eqref{eq:cVdP}. For a comprehensive discussion of the folded singularities in our model \eqref{eq:cVdP}, see \ref{app:gspt}, where we employ the geometric singularity perturbation theory \cite{Fenichel1979}.

\subsection{Strong symmetry breaking: from double-loop oscillations to single-loop endocycles}\label{sec:symmetrybreaking}

\begin{figure}[htp!]
    \centering   
    \includegraphics[width=\textwidth]{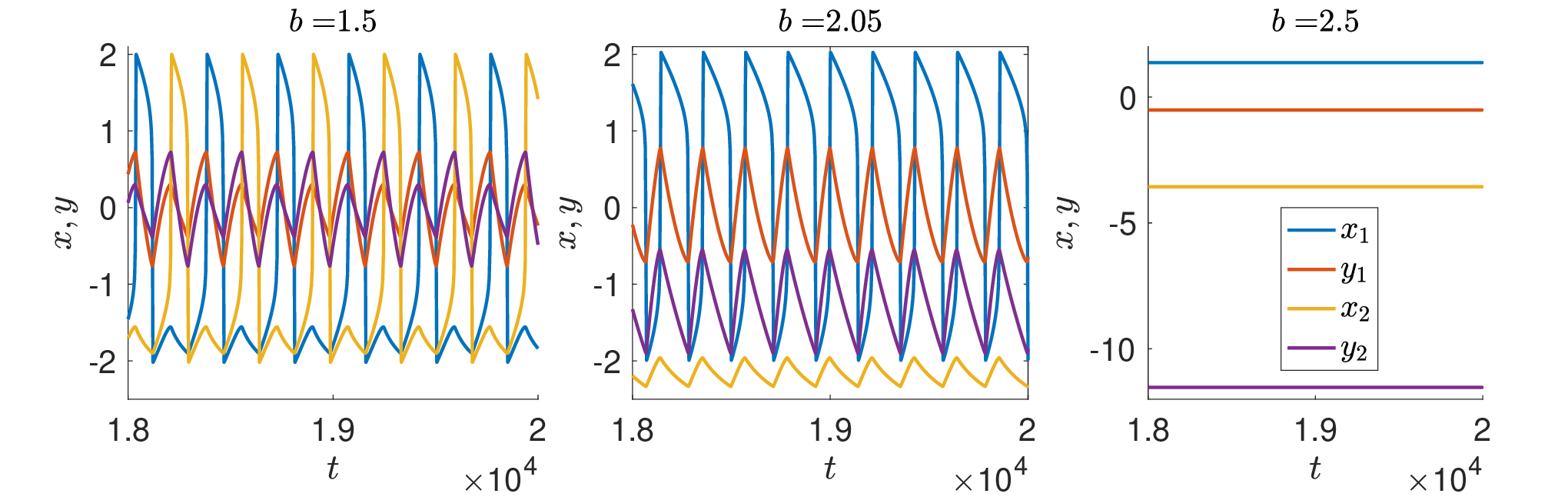}
    \caption{Increasing the symmetric coupling strength $b$ breaks the symmetry in the dynamics and leads to the transition from the double-loop oscillations (A,B,C) to the single-loop oscillations }
    \label{fig:transition_b}
\end{figure}

In \cite{Dragoi2024cellcycle}, the authors varied one of the two coupling strengths to place one oscillator in the latched state while allowing the other to oscillate. Unsurprisingly, this results in the transition from double-loop oscillations to single-loop endocycles. 
In this paper, we demonstrate that a similar transition can be induced in \eqref{eq:cVdP}, but with a key distinction: we symmetrically vary the coupling strength $b$, increasing the inhibition to both oscillators equally (see Figure \ref{fig:transition_b}, left and middle panels). This reveals a novel and interesting mechanism underlying the strong symmetry breaking for the double-oscillation/endocycle transition, as will be explained in the next subsection using bifurcation analysis. For sufficiently large values of inhibition $b$, the system transitions to an asymmetric steady state, characterized by one oscillator having a higher $x$ value than the other (Figure \ref{fig:transition_b}, right).

\subsubsection{Transitions as a function of symmetric coupling strength}\label{sec:onepara}



In this section, we explore the mathematical mechanism underlying the transitions between double-loop oscillations that represent normal cell functioning and  single-loop oscillations that represent endocycles. Our numerical bifurcation analysis suggests a homoclinic bifurcation within this $\mathbb{Z}_2$-symmetric model is responsible for the transition.  

\begin{figure}[t!]
    \centering   
   \includegraphics[width=\textwidth]{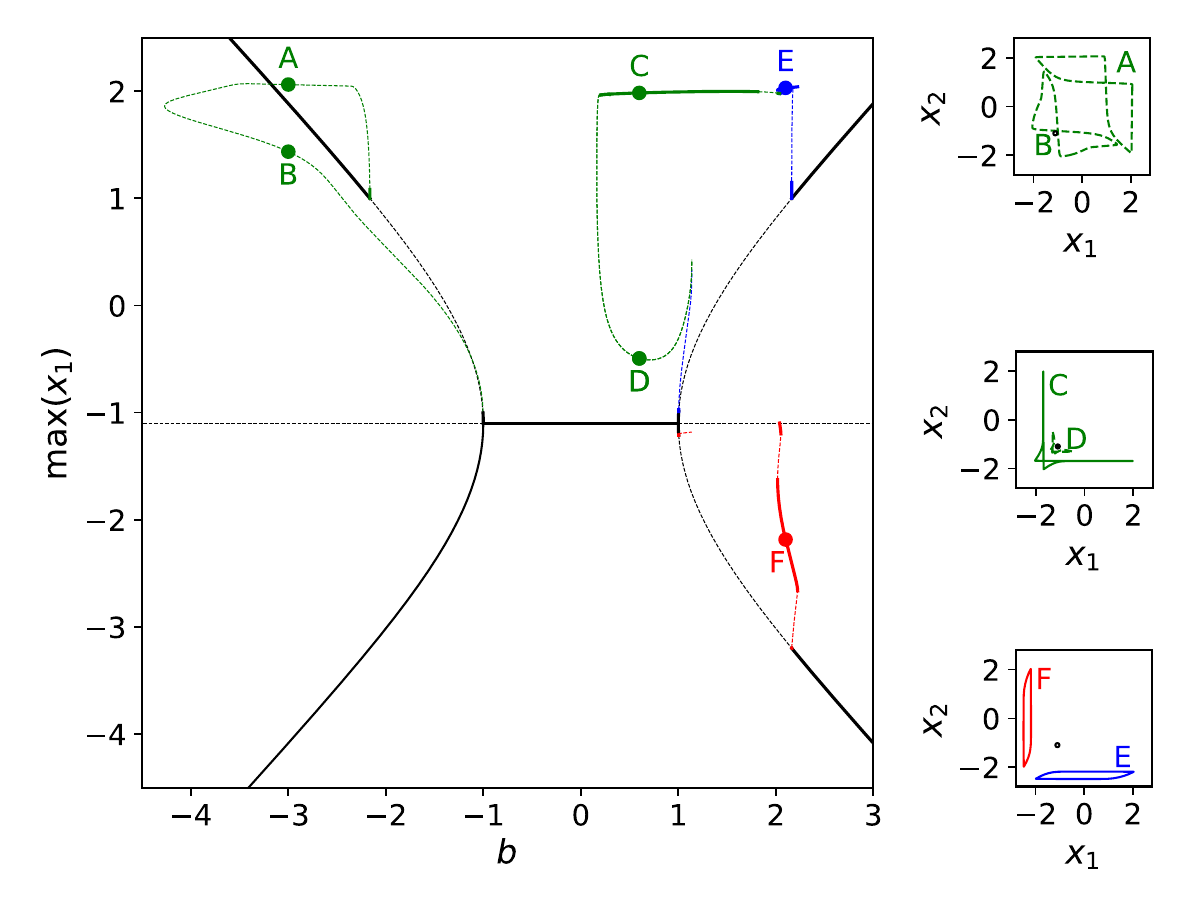}\\
    \caption{Bifurcation diagram of \eqref{eq:cVdP} with respect to the coupling strength $b$, consisting of steady-state equilibria (black), double-loop periodic orbits (green) and single-loop periodic orbits (red and blue). Solid (resp., dashed) curves denote stable (resp., unstable) branches. The representatives solutions corresponding to parameters labeled A, B, C, D, E and F from the left panel are projected onto $(x_1,x_2)$-plane in the right panel.  
    }
    \label{fig:bifb}
\end{figure}

Figure~\ref{fig:bifb} shows a partial numerically-computed bifurcation diagram of the system~\eqref{eq:cVdP} as a function of the coupling parameter $b$ with other parameters fixed at Table \ref{tab:parameters},
along with projections of trajectories in the $(x_1,x_2)$-plane for several Representative solutions. Steady-state equilibria are shown in black, double-loop periodic orbits in which the oscillators are out of phase are shown in green, and single-loop periodic orbits where one oscillator has a significantly larger amplitude than the other (i.e., endocycles) are shown in blue and red.  For periodic orbits, the maximum value of $x_1$ is shown and, in all cases, stability (instability) is indicated by thick solid (thin dotted) lines.  All computations are carried out using the numerical continuation package AUTO~\cite{doedel1981auto}.

As the magnitude of the coupling strength increases from $b=0$ to $b=\pm 1$, we see a loss of stability and the emergence of additional pairs of stable steady-states through pitchfork bifurcations.  These emerging states still satisfy $x_1=x_2$ at the $b=-1$ pitchfork, while the $b=+1$ pitchfork breaks this symmetry.  We focus on our attention on positive coupling parameter $b>0$ (i.e., mutual inhibition) as the relevant case for capturing cell cycle dynamics described in Figures~\ref{fig:biomodel} and \ref{fig:soln}.    While $x_1\neq x_2$ along the branches emerging from $b=+1$, these branches map onto each other under exchange of $x_1$ and $x_2$.  Supercritical Hopf bifurcations occurs on both asymetric steady-state branches just after the pitchfork bifurcations at $b=1.0041$, after which the steady-states lose stability and a symmetric pair of single-loop periodic orbit branches are born (see the blue and red curves near $b=1$ in Figure \ref{fig:bifb}). These single-loop periodic orbits  quickly lose stability at $b=1.0042$ through a torus bifurcation and the period of the orbits begins to diverge as $b$ approaches a homoclinic (HC) bifurcation at $b=1.1379$ associated with asymmetric steady-states. 
The steady-states regain stability through another supercritical Hopf bifurcation that occurs near $b=2.164$, at which another symmetric pair of single-loop periodic orbit branches indicated by blue and red curves are formed. These single-loop periodic orbits lose and gain stability through saddle-nodes of periodic orbits and representative stable orbits labeled E and F are shown in the lower right panel. The period eventually begins to diverge as $b$ approaches another homoclinic bifurcation associated the symmetric equilibrium at $b= 2.0392$.   

\begin{figure}[t!]
    \centering   
    \includegraphics[width=\textwidth]{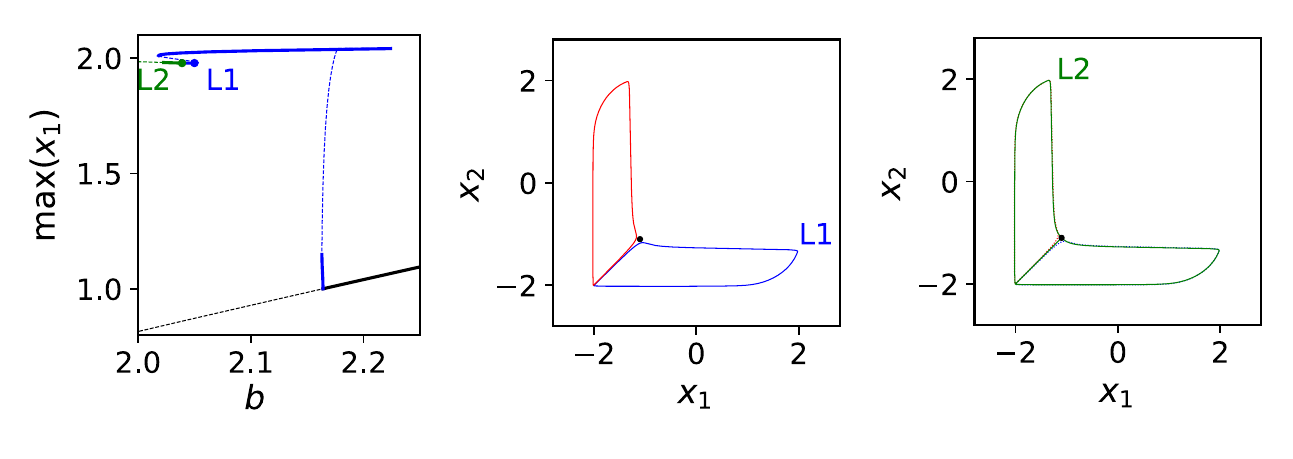}
    \caption{Enlarged view of the bifurcation diagram from Figure \ref{fig:bifb} near the Homoclinic bifurcation along the upper periodic orbit branch, with representative solutions for parameters on either side of the homoclinic bifurcation (L1: $b=2.05$ and L2: $c=2.039$).
    }
    \label{fig:bifb_zoom}
\end{figure}

Initializing continuation using the solution shown in Figure~\ref{fig:soln}, we also compute a double-loop branch of periodic orbits with representative solutions labeled C and D shown in the middle panel on the right. The period of the orbits diverge on either side of the branch as the double-loop branch approaches the single-loop branch at a homoclinic bifurcation. The left panel of Figure~\ref{fig:bifb_zoom} zooms in near the location where the double-loop branch meets the upper single-loop branch. As we approach the homoclinic point from the right, the two single-loop orbits on the red and blue branches of Figure~\ref{fig:bifb} approach the unstable symmetric steady state. Projections of these periodic orbits and the steady states at the point labeled L1 are shown in middle panel of Figure~\ref{fig:bifb_zoom}. On the left of this homoclinic point, a single double-loop orbit is formed where the two single-loop orbits used to be. A projection of this orbit at the point labeled L2, along with the steady state, is shown in the right panel.  The mechanism underlying this transition is explored further in section~\ref{sec:HCButterfly}.  

\subsubsection{Homoclinic bifurcation and strong symmetry breaking}\label{sec:HCButterfly}

\begin{figure}[!t]
\begin{center}
\begin{minipage}{.75\textwidth}
\begin{tabular}
{@{}p{0.5\linewidth}@{\quad}p{0.5\linewidth}@{}}
\subfigimg[width=\linewidth]{\bfseries{\small{(A)}}}{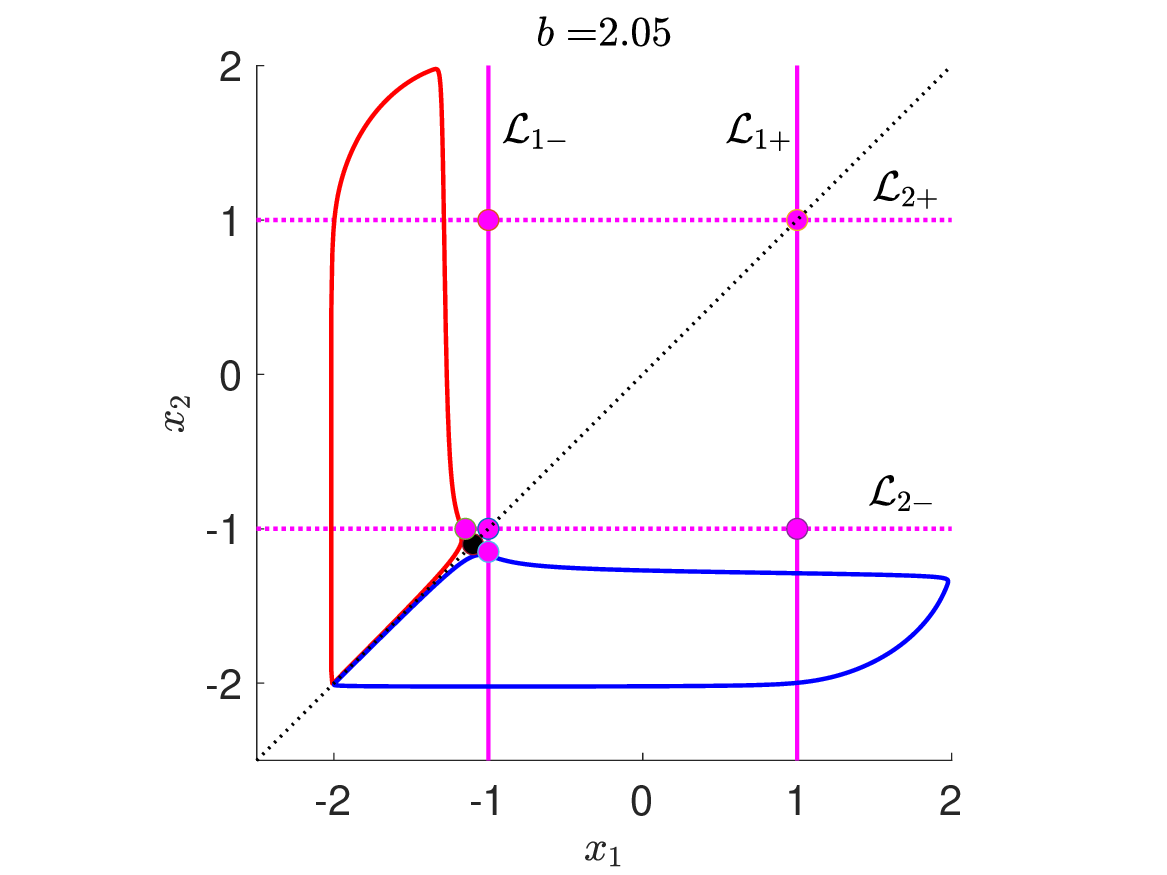} & \subfigimg[width=\linewidth]{\bfseries{\small{(B)}}}{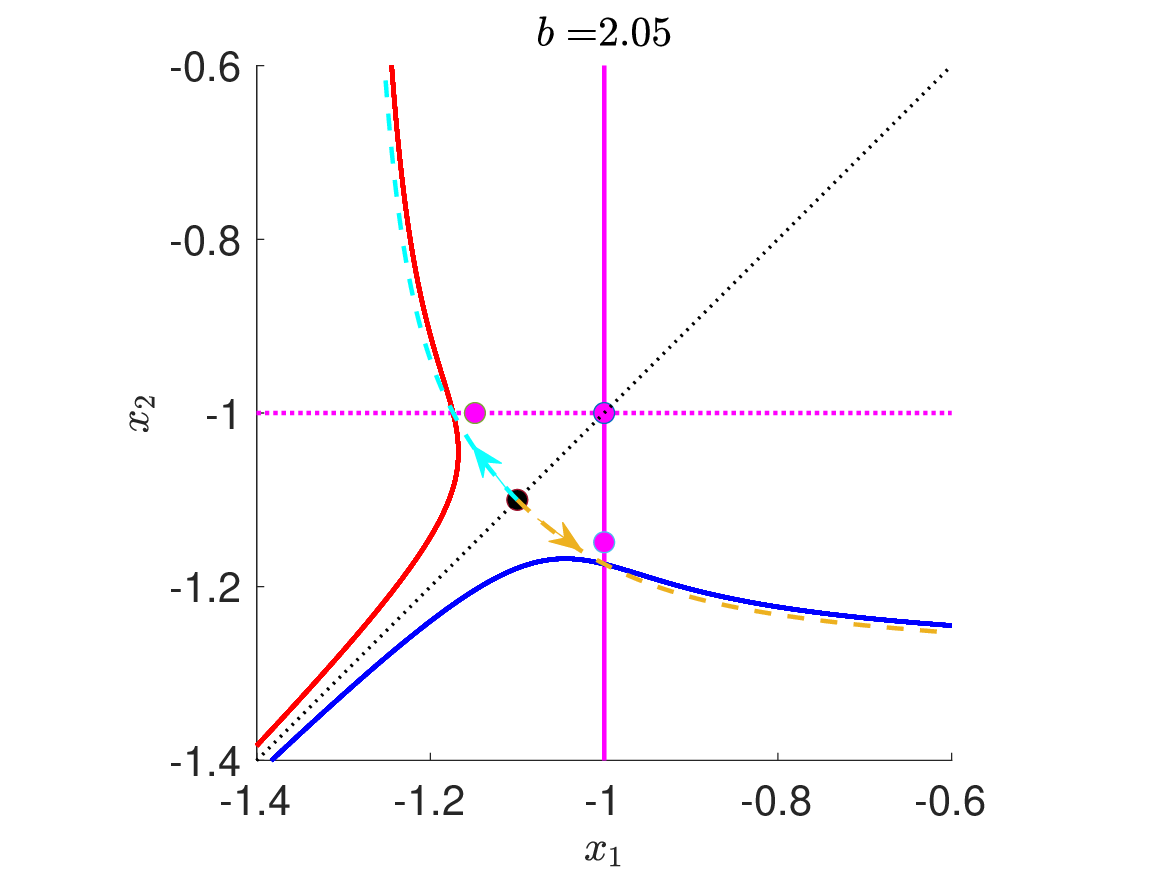}\\
\subfigimg[width=\linewidth]{\bfseries{\small{(C)}}}{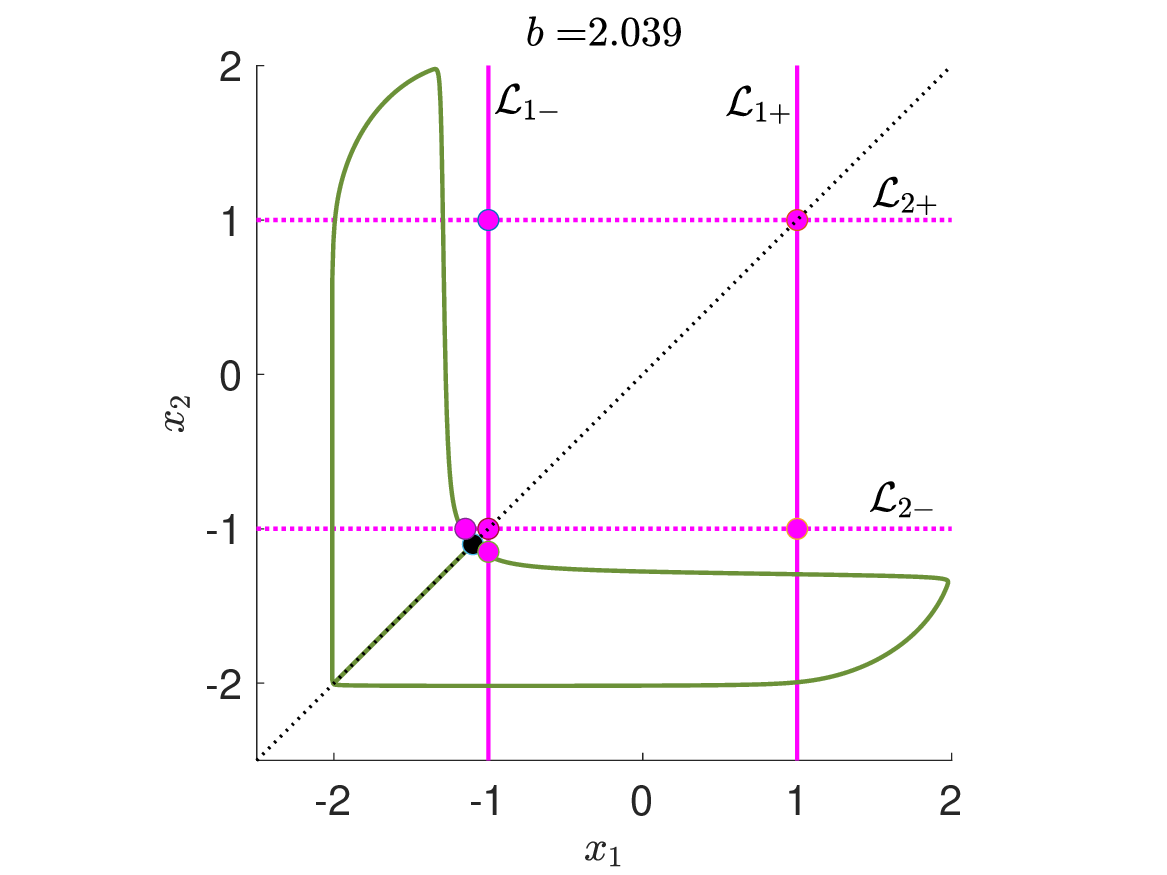}&
\subfigimg[width=\linewidth]{\bfseries{\small{(D)}}}{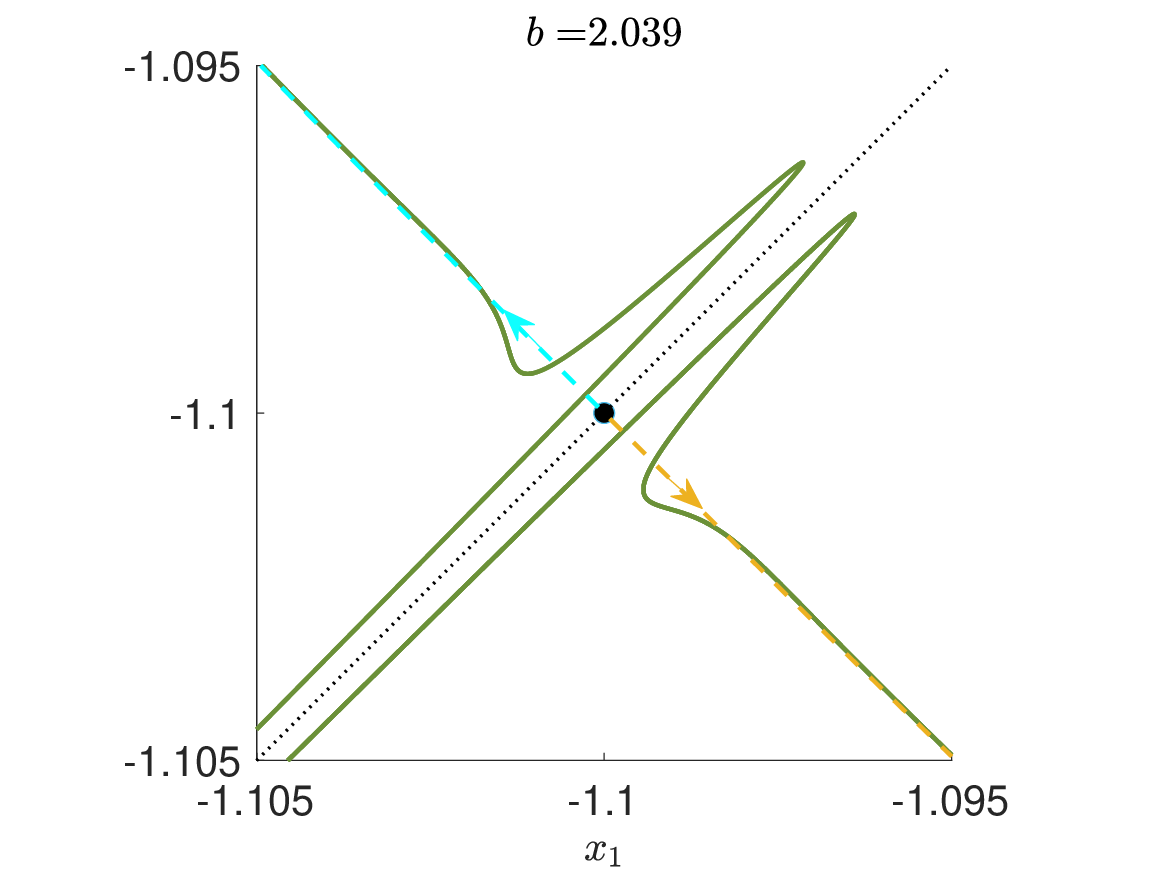}\\
\multicolumn{2}{c}{\subfigimg[width=\linewidth]{\bfseries{\small{(E)}}}{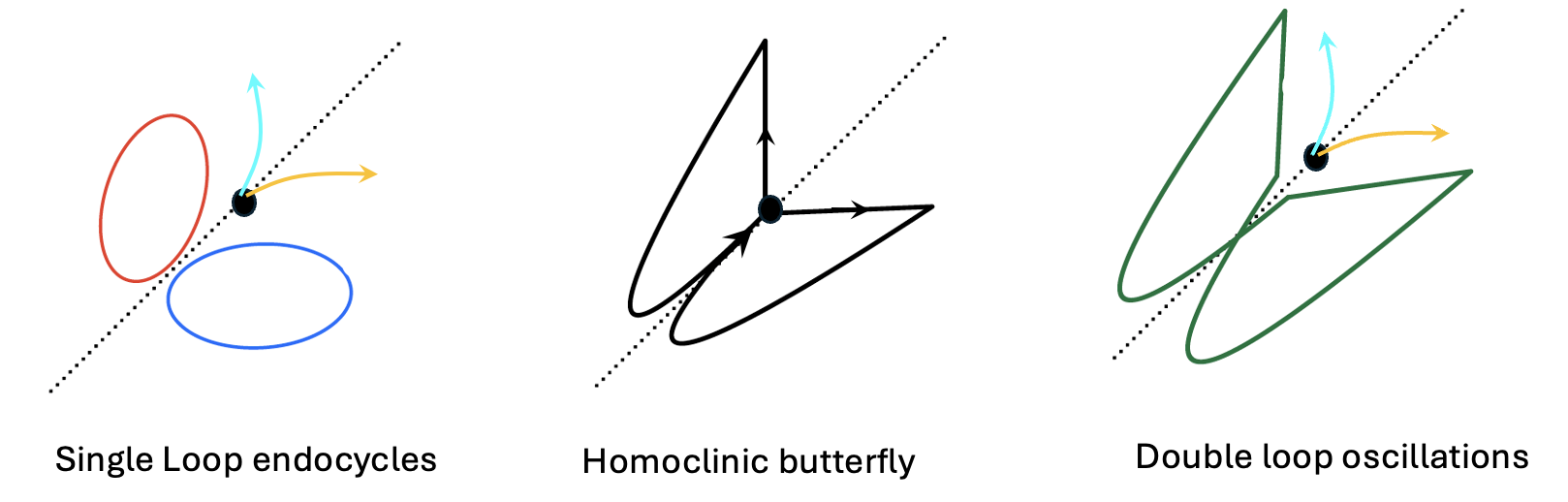}}
\end{tabular}
\end{minipage}
\end{center}
\caption{Increasing the symmetric coupling strength $b$ across the homoclinic bifurcation breaks the symmetry in the double-loop periodic orbit. (A,B) The single-loop oscillations at $b=2.05$ corresponding to L1 in Figure \ref{fig:bifb_zoom}. (C,D): The double-loop oscillations at $b=2.039$ corresponding to L2 in Figure \ref{fig:bifb_zoom}. The red and blue trajectories are two different single-loop oscillations from different initial conditions. (B) and (D) show enlarged views of (A) and (C) near the symmetric equilibrium $E_0$ (black circle). The cyan and yellow dashed trajectories represent the two branches of the unstable manifold of $E_0$. For both $b=2.039$ and $b=2.05$, the two asymmetric folded singularities (magenta circles) near $E_0$ are stable folded foci and the one at $(-1, -1)$ is a folded saddle. Other coding coding and symbols have the same meanings as in Figure \ref{fig:soln}C. (E): Schematic representations of single-loop endocycles, homoclinic butterfly, and double-loop oscillations. Hand-drawn trajectories illustrate key details near the symmetric axis (black dashed line), which are difficult to visualize in the actual solutions in (A-D) due to scale limitations. 
}
\label{fig:transit-to-endo}
\end{figure} 

Periodic orbit solutions labeled L1 and L2 from Figure \ref{fig:bifb_zoom} are shown again in Figure \ref{fig:transit-to-endo}A and C, along with the folds $\mathcal{L}_{i\pm}$ (magenta lines) and all folded singularities (magenta circles). Panels (B) and (D) provide enlarged views of the dynamics near the symmetric equilibrium $E_0$, which is of saddle type (black circle). The cyan and yellow trajectories are two branches of the unstable manifold of $E_0$, denoted as $W^U(E_0)$. For both $b=2.05$ and $b=2.039$, the two asymmetric folded singularities near $E_0$ are folded foci instead of folded saddles as in the previous case with $b=0.3$ (Figure \ref{fig:soln}C). These singularities do not appear to play a significant role in the dynamics or contribute to symmetry breaking, in contrast to what is observed in \cite{awal2023symmetry,awal2024symmetry}. Instead, the mechanism driving the transition from double-loop to single-loop oscillations (i.e., strong symmetry breaking) in our system involves a double homoclinic bifurcation, also referred to as homoclinic butterfly (see e.g., \cite{shilnikovL2003HC}). 

We illustrate the homoclinic mechanism underlying the transition using schematic trajectories in Figure \ref{fig:transit-to-endo}E. The left and right panels in (E) correspond to the solutions from the top two rows, with consistent color coding.  In the endocycle scenario, the upper branch of $W^U(E_0)$, depicted by the cyan trajectory, converges toward the upper single-loop periodic orbit (red curve), whereas the lower branch shown by the yellow trajectory converges to the lower single-loop periodic orbit solution (blue curve). 
At the homoclinic bifurcation, the two branches of $W^U(E_0)$ are double asymptotic to $E_0$. Specifically, they approach $E_0$ tangentially to each other along the symmetric axis (see Figure \ref{fig:transit-to-endo}E, middle panel). After the bifurcation, the two branches of $W^U(E_0)$ twist and intersect along the symmetric axis when projected onto $(x_1,x_2)$, forming double-loop oscillations, as illustrated in the right panel of (E). 

\subsection{Transitions as a function of symmetric coupling strength and latching parameter}\label{sec:twopara}

We now analyze the dependence of the sequence of transitions between single-loop oscillations, double-loop oscillations and steady state equilibria through a combination of simulation and numerically-computed two-parameter bifurcation diagrams.  We focus on understanding how the bifurcations, as a function of coupling parameter $b$, depend on the latching parameter $a$.

\begin{figure}[t!]
    \centering   
    \includegraphics[width=0.9\textwidth]{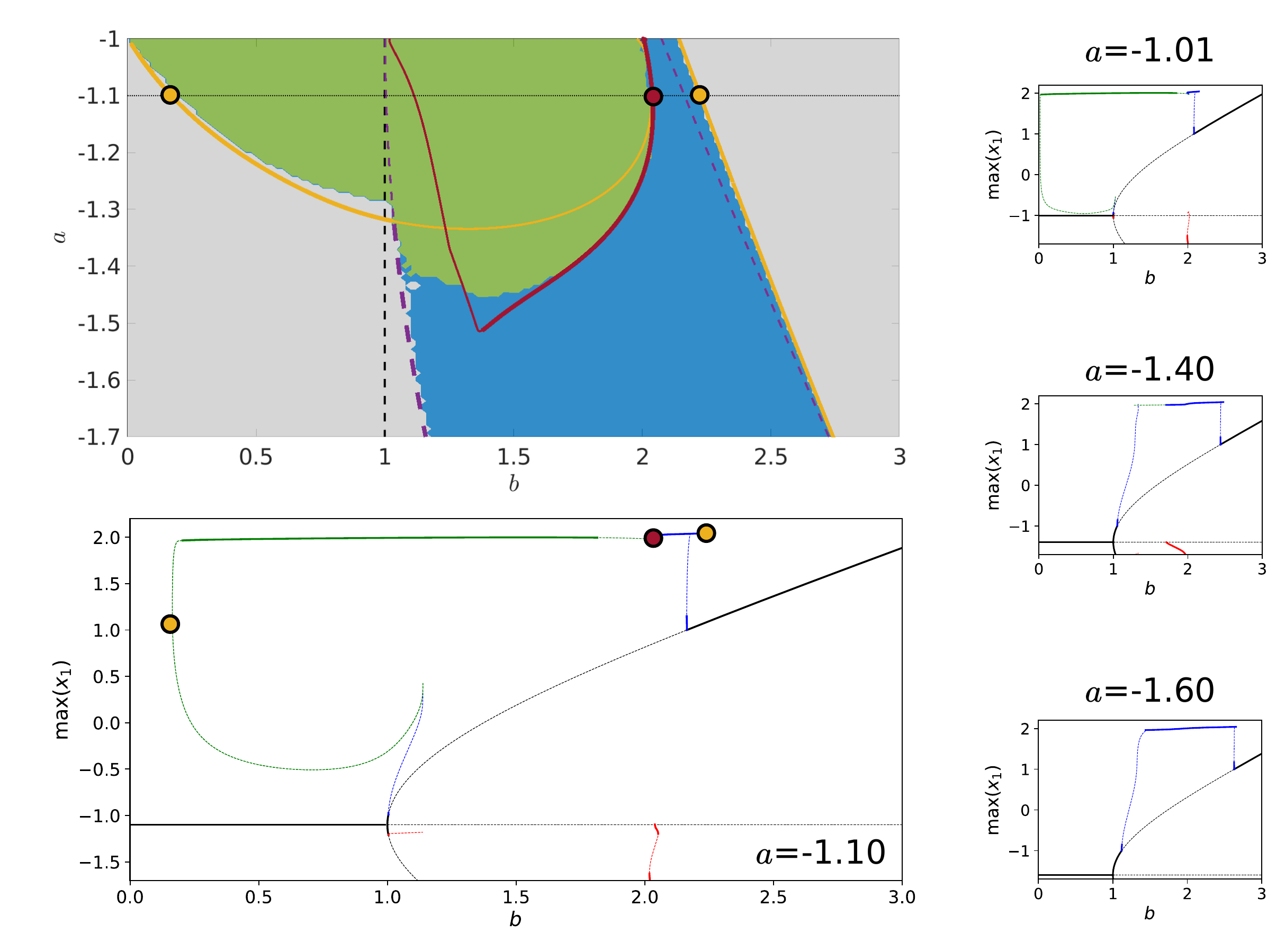}
    \caption{Two-parameter bifurcation diagram of \eqref{eq:cVdP} with respect to $(a,b)$. Curves of left SNPO of the double-loop periodic orbits (yellow solid), Hopf (purple dashed), homoclinic (red solid) and the right-most SNPO of the single-loop periodic orbits divide the space into different solution patterns observed in simulation. Regions with steady states are light gray, single-loop in blue and double-loop in green. Thick lines are expected to align with the observed boundary between different states. 
    Representative one-parameter bifurcation diagrams with respect to $b$ for $a=-1.6, -1.4, -1.1, -1.01$ are also shown. Color coding of the one-parameter bifurcation diagram is the same as Figure \ref{fig:bifb}. }
    \label{fig:pattern-ab}
\end{figure}

We conduct numerical simulations of the system on a grid of parameter values with $-1.7<a<-1$, $0<b<3$ and all other parameters given in Table~\ref{tab:parameters}.   We restrict ourselves to $a<-1$  so that the steady-state equilibrium with $x_1=x_2=a$ remains stable in the uncoupled case of $b=0$.  Each simulation is initialized with a large perturbation in the $x_2$ coordinate of steady-state equilibrium, and the final state of the system is classified as either steady-state, single-loop oscillation or double-loop oscillation after $10^4$ units of time.   We distinguish single loops from double loops by checking if just one or both $x_i$ undergo large amplitude oscillations as illustrated in Figure~\ref{fig:transition_b}.  We also check relative phase of the oscillators to ensure that they are not synchronous in the double-loop region.  The results are indicated by the shading in the upper left panel of Figure~\ref{fig:pattern-ab} with steady states in light gray, single loop oscillations in blue and double loop oscillations in green.

A one-parameter bifurcation diagram in $b$ is shown for $a=-1.1$ in the lower left panel of Figure~\ref{fig:pattern-ab}, which corresponds to the upper right section of Figure~\ref{fig:bifb}.  Additional one-parameter bifurcation diagrams are shown on the right with at $a=-1.01$, $-1.4$ and $-1.6$ in order from top to bottom.  We additionally compute the location of the following bifurcations in the $(a,b)$ parameter plane:
\begin{itemize}
    \item Black dashed line: pitchfork bifurcation
    \item Purple dashed lines: Hopf bifurcations (HB)
    \item Orange solid lines: Left saddle-node of double-loop periodic orbits (SNPO) and right-most SNPO of the single-loop periodic orbits 
    \item Red solid line: Homoclinic bifurcations 
\end{itemize}
We make the sections of these two-parameter bifurcation curves thick where we expect them to agree with the transitions between states obtained from numerical simulation.  Note that we approximate the Homoclinic bifurcation curve by numerical continuation of a periodic orbit near the bifurcation as a function of $a$ and $b$  with the period fixed to a sufficiently long value.

At a fixed value of $a=-1.1$, we see the appearance of double-loop oscillations in simulations as $b$ increases past the left SNPO (orange line). There is a region of bistability between the symmetric steady state and the double-loop periodic orbits, which ends when this steady state loses stability through a pitchfork bifurcation (dashed black line).  There is a small section of bistability between the emerging asymmetric steady states and the double-loop periodic orbits between this pitchfork and the lower Hopf bifurcation (dashed purple line).  The periodic orbits near the lower homoclinic bifurcation are unstable and therefore do not mark transitions in the states observed in simulation.  The upper homoclinic orbit, on the other hand, aligns with the transition between double-loop and single-loop oscillations seen in the simulation data. The single-loop oscillations persist until the left-most SNPO (solid orange line) on the single-loop branch.  There is a small region of bistability between the stable asymmetric steady states and the single-loop periodic orbits before only the asymmetric stable steady-states remain.  

For $a$ closer to $a=-1$, where the uncoupled oscillators become unlatched and undergo Hopf bifurcations, the region where double-loop oscillations can be found stretches down to near $b=0$.  This is illustrated by the one-parameter bifurcation diagram with $a=-1.01$ on the right. As $a$ decreases away from the Hopf bifurcation at $a=-1$ and the system becomes less excitable, we see a decrease in the range of $b$ values that support double-loop oscillations (green region).  We also see an increase in the region where single-loop oscillations (blue region) can be found. As $a$ decreases to $a=-1.4$, the one-parameter bifurcation on the right indicates that there are no stable double-loop periodic orbits.  While we do see double-loop oscillations in the simulations, very long runs indicate that they are not periodic and undergo occasional disruptions to the nearly-periodic behavior.  Eventually, the double-loop oscillations disappear from simulations for $a$ value below the fold of homoclinic bifurcation curve (red solid) and only single loops remain. As the one-parameter bifurcation at $a=-1.6$ illustrates, the two Hopf bifurcations along the asymmetric steady states are now connected by a single branch of single-loop periodic orbits.   



\section{Discussion}\label{sec:discuss}
In this work, we have proposed identical coupled van der Pol oscillators as a conceptual model for the latching mechanism in cell cycling as described by~\cite{novak2022mitotic}.  Each of the two oscillators is biased to just below their Hopf bifurcation so that there is only a globally attracting steady state when uncoupled.  It is in this sense we consider the regulatory networks ``latched."  These oscillators, however, sit in an excitable regime so that a sufficiently large perturbation can induce a single large-amplitude oscillation before returning to the latched state. \RED{A similar observation—that coupling two excitable systems can generate large-amplitude, out-of-phase oscillations—dates back to early work on echo waves~\cite{krinskii1972investigation, tyson1979oscillations}, which demonstrates that mutual excitation between such systems can cause a pulse of activity to `echo' back and forth. Oscillations of this type have been observed in various contexts, including heart tissue, neural tissue, and chemical reactions \cite{winfree1972spiral,belakhovski1965several,gul1972mechanism,winfree1974rotating,winfree1980geometry}.} 
Motivated by cell cycle biology,  we consider inhibitory coupling so that one oscillator suppresses the other during the positive part of its oscillation and excites the other to initiate a large oscillation during the negative part of its oscillation.  In this way, mutual inhibition, parameterized by coupling strength $b$, leads to instability of the latched steady state and produces sustained double-loop oscillations of the type shown in Figure~\ref{fig:soln}, which represent normal cell cycling in \cite{novak2022mitotic} (Figure~\ref{fig:biomodel}).

We identify, through a combination of simulation and numerical continuation, a transition to single-loop endocycles in which only one of the oscillators undergoes sustained large-amplitude oscillations (c.f., Figure~\ref{fig:transition_b}). This strong symmetry-breaking transition occurs as the symmetric inhibitory coupling strength $b$ becomes even stronger.  We note that endocycles are possible in this model even when both oscillators are in the latched state (e.g. $a<-1$).  The underlying mathematical mechanism responsible for the transition to endocycles, as schematically illustrated in Figure~\ref{fig:transit-to-endo}(E), is a homoclinic bifurcation. \RED{For $a=-1.1$, the single-loop oscillations occur only within a narrow range of the coupling strength $b$ (see Figures~\ref{fig:bifb} and~\ref{fig:bifb_zoom}, blue periodic orbit branch). This poses a challenge in associating them with robust biological events, such as single-loop endoreplication \cite{edgar2001}. However, this range (blue region in Figure~\ref{fig:pattern-ab}) broadens with decreasing values of the latching parameter $a$, which controls the transition between latching and swinging within the individual oscillators. Conversely, double-loop oscillations (green region in Figure~\ref{fig:pattern-ab}) become less common as the oscillators are pushed to a more strongly latched state, corresponding mathematically to decreasing the excitability of the oscillators.}
\RED{We note that in our symmetrically coupled model, both the single-loop G1/S $(x_1,y_1)$-oscillations, which can be associated with endoreplication cycles \cite{edgar2001,novak2022mitotic}, and single-loop M/G1 $(x_2,y_2)$ oscillations, which can be associated with Cdc14 endocycles \cite{lu2010,manzoni2010,novak2022mitotic}, emerge from the same underlying mechanism due to symmetry. However, this symmetry may not align with biological reality, where these processes are likely driven by distinct mechanisms. Furthermore, while endoreplication represents a robust biological event, 
Cdc14 endocycles are rare pathological states. This limitation highlights the potential of exploring asymmetric coupled oscillator models to better capture the biological specificity of cell cycle regulation. }

We have been careful to make a distinction between sustained oscillations and periodic orbits within the range of dynamics predicted by this model.  In particular, we observe double-loop oscillations that are not periodic in simulations with parameter values within the green region of Figure~\ref{fig:pattern-ab}.  These oscillations are nearly periodic in the sense that they alternate in a double-loop fashion for tens to hundreds of cycles or more between disruptions.  From a biological standpoint, this may not be a major issue as it predicts normal cell cycling for a vast majority of the time. More problematic for the biology are cases in which a sequences of large oscillations occur in one oscillator before switching to a sequence of large oscillations in the other.  These type of solutions were observed occasionally in simulations, though were restricted to narrow areas of the parameter space where   just to the left of the homoclinic bifurcation responsible for the transition between double-loop oscillations and single-loop endocycles.

Our bifurcation analysis highlights an important role of the latching parameter $a$ in controlling the cell division cycle (see Figure \ref{fig:pattern-ab}). We focus on the region where $a<-1$, which biases the oscillators towards a latched state. For $a$ values close to $-1$, a bistable region forms between the left SNPO curve and the pitchfork bifurcation. In this scenario, increasing the symmetric coupling $b$ induces a sequence of transition from a symmetric steady state to double-loop oscillations, then to single-loop endocycles, and eventually to an asymmetric steady state.  As $a$ decreases further (e.g., $a=-1.4$), making the system more latched in isolation, the transition to double-loop oscillations occurs at a Hopf bifurcation rather than an SNPO, causing the bistable region to vanish. With $a$ falling below the fold point of the homoclinic bifurcation curve, the double-loop oscillations associated with normal cell cycling are lost. This analysis suggests that a model of the cell cycle's latching mechanism must carefully balance the degree of latching to maintain normal function; excessive latching disrupts the cell cycle progression, while too little can excessively extend the double-loop oscillation phase, \RED{resulting in a limited region of endoreplication, as discussed above}. These thresholds for state transitions and bistability could be experimentally tested by manipulating the latching parameter, coupling strength, and perturbing cell cycle dynamics. 

Recently, coupled oscillators were also used in \cite{Dragoi2024cellcycle} to capture the latching dynamics in a complex cell cycle model in human cells introduced in \cite{Dragoi2024model}.  \RED{Their system is closely tied to the known biochemistry of eukaryotic cell cycle regulation, making it much easier to associate their simulations with experimental observations. Our model does not aim to capture detailed biological processes but instead seeks to qualitatively reproduces observed behaviors and examine the mechanisms underlying transitions between dynamics described in \cite{novak2022mitotic}. Our choice of parameters in the coupled van der Pol model~\eqref{eq:cVdP} puts the uncoupled oscillators in an excitable regime in contrast to ~\cite{Dragoi2024cellcycle}.  This difference leads to a mathematically distinct mechanism for the transition from regular cell cycling to endocycles. In particular, we have identified a homoclinic bifurcation as a mechanism for the strong symmetry breaking responsible for this transition.  Exploring more biologically grounded cell cycle models (e.g. \cite{novak2022mitotic,Dragoi2024cellcycle,Dragoi2024model}) for evidence of a similar homoclinic bifurcation mechanism represents an interesting future direction,  and may lead to deeper insights into the potential relevance of this mechanism to cell cycle regulation.  }

Strong symmetry breaking analogous to double-loop oscillation/endocycle transition have also recently been observed in other coupled oscillator systems \cite{awal2019symmetry,awal2020symmetry,awal2021symmetry,awal2023symmetry,awal2024symmetry,epstein2024}. References~\cite{awal2023symmetry,awal2024symmetry,epstein2024} show that the strong symmetry-breaking mechanism is due to folded singularities near the symmetric equilibrium. Specifically, they demonstrate that a folded node located off the symmetry axis is the key mechanisms in the coupled Lengyel-Epstein oscillators~\cite{awal2024symmetry}.
However, in our system, folded singularities do not play a significant role in the dynamics. In fact, folded nodes are not present within the considered parameter regime. See \ref{app:gspt} for more details. Similar to the coupled oscillators in \cite{Dragoi2024model}, the symmetric equilibria of the coupled oscillators in ~\cite{awal2023symmetry,awal2024symmetry} are also unstable. This contrasts with our system, where the symmetric equilibrium is stable, leading to differences in the observed dynamics. Additionally, the type of coupling in our model is different: we employ nonlinear fast-to-slow coupling, whereas theirs involve linear fast-to-fast or slow-to-slow couplings. 


Our focus in this work has been developing the the coupled van der Pol model as a conceptual framework for understanding the latching mechanism for cell cycle regulation.  There are, however, a very rich range of dynamical behaviors predicted by this simple model that we be explored in future work. In addition to double-loop periodic orbits in which the oscillators are out of phase, simulations produce nearly periodic trajectories in addition to more complicated dynamics.  The fold point on the homoclinic bifurcation curve of Figure~\ref{fig:pattern-ab} may serve as an organizing center for the various dynamics related to double-loop oscillations in the model.  From a biological perspective, it will be important to also explore dynamics in the case that the oscillators are no longer identical.    







\section*{Acknowledgments}
This work was initiated through discussions during the 
2023 conference on \textit{Dynamical Systems in the Life Sciences} at The Ohio State University, Columbus, OH.   
This work was supported in part by NIH/NIDA R01DA057767.

\section*{Conflict of interest}

The authors declare there is no conflict of interest.

\newpage
\appendix

\section{Geometric Singular Perturbation Theory and Folded Singularity}\label{app:gspt}

In this section, we apply the Geometric Singular Perturbation Theory (GSPT) \cite{Fenichel1979,Nan2015} to the coupled fast-slow oscillators \eqref{eq:cVdP} by treating $\varepsilon$ as the singular perturbation parameter. Our goal is to identify and classify all the folded singularities \cite{awal2024symmetry,PW24}. 

Taking $\varepsilon\to 0$ in \eqref{eq:cVdP} yields the two-dimensional (2D) \emph{fast layer problem}
\begin{subequations}
\begin{align}
    \dot{x_1} &=y_1 + \left(1 -\frac{x_1^2}{3}\right)x_1=: f_1 \\
    \dot{x_2} &=y_2 + \left(1 -\frac{x_2^2}{3}\right)x_2=: f_2
\end{align}
\end{subequations}
where $y_1$ and $y_2$ are constant. 

Equilibria of the fast subsystem gives the critical manifold:
\begin{align}\label{eq:ms}
    \mathcal{M_S} = \left\{y_1 = (-1 +\frac{x_1^2}{3})x_1 ,\,
    y_2=  (-1 +\frac{x_2^2}{3})x_2   \right\}
\end{align}

Stability of the $\mathcal{M_S}$ is determined by the eigenvalues of the Jacobian of the fast layer problem evaluated along $\mathcal{M_S}$,

\begin{equation}\label{eq:J-matrix}
J=
\begin{pmatrix}
f_{1,x_1} & f_{1,x_2}  \\
f_{2,x_1} & f_{2,x_2} \\
\end{pmatrix}
=
\begin{pmatrix}
1-x_1^2 & 0  \\
0 & 1-x_2^2 \\
\end{pmatrix}
\end{equation}
It follows that the eigenvalues of $J$ are real, and given by $1-x_1^2$ and $1-x_2^2$. Thus, $\mathcal{M_S}$ contains only stable and unstable nodes, saddles and saddle-nodes. Boundaries of different equilibria are given by the subset of $\mathcal{M_S}$ on which $\det(J)=0$, which is called the fold curve
\begin{subequations}\label{eq:fold}
\begin{align}
    \mathcal{L_S} &= \left\{\det(J)=0 \right\} \\&= \left\{ (1-x_1^2)(1-x_2^2)=0  \right\} \\& = \left\{  x_1=\pm 1, x_2 = \pm 1 \right\}\subset \mathcal{M_S}.
\end{align}
\end{subequations}
Together, there are four fold lines $\mathcal{L}_{1\pm}$ and $\mathcal{L}_{2\pm}$, where $\mathcal{L}_{1-}=\{x_1=-1\}$ and $\mathcal{L}_{1+}=\{x_1=+1\}$ are the left and right folds of $x_1$-nullcline while $x_2$ can be any value, $\mathcal{L}_{2-}=\{x_2=-1\}$ and $\mathcal{L}_{2+}=\{x_2=+1\}$ are the left and right folds of $x_2$-nullcline while $x_1$ can be any value.

Introducing a slow time $t_s = \varepsilon t$ yields an equivalent descrition of dynamics of \eqref{eq:cVdP}:
\begin{subequations}\label{eq:slowcVdP}
\begin{align}
    \varepsilon_1 x'_1 &=y_1 + \left(1 -\frac{x_1^2}{3}\right)x_1\\
    y'_1 &=(a_1-x_1 + b_{1} \tanh(k_2(a_2-x_2))) \\ 
   \varepsilon_2 x'_2 &=y_2 + \left(1 -\frac{x_2^2}{3}\right)x_2\\
    y'_2 &= (a_2-x_2 + b_{2}\tanh(k_1(a_1-x_1))) 
 \end{align}
\end{subequations}
where $'$ denote the derivative with respect to the slow time $t_s$. 

Taking the same singular limit $\varepsilon\to 0$ in the slow system \eqref{eq:slowcVdP} yields the 2D \emph{slow reduced problem}, a system that describes the dynamics of $y_1,y_2$ along $\mathcal{M_S}$,
\begin{subequations}\label{eq:slow-sub}
\begin{align}
    y_1'&=a_1-x_1 +b_{12} \tanh(k_1(a_2-x_2)) =: g_1(x_1,x_2)\\
    y_2' &=a_2-x_2 +  b_{21} \tanh(k_1(a_1-x_1)) =: g_2(x_1,x_2)
\end{align}
\end{subequations}
where $(x_1,y_1,x_2,y_2)\in \mathcal{M_S}$.

To investigate canard dynamics due to folded singularity, we project the slow reduced problem \eqref{eq:slow-sub} onto the fast-slow subspace to obtain a complete description of the dynamics along $\mathcal{M_S}$. To this end, we differentiate the graph representation of $\mathcal{M_S}$ in \eqref{eq:ms} to obtain
\begin{subequations}\label{eq:reduced-slow-sub}
\begin{align}
   (x_1^2-1) x_1' &=g_1\\
   (x_2^2-1) x_2' &=g_2\\
    y_1'&=g_1\\
    y_2' &=g_2 
\end{align}
\end{subequations}
Note that the reduced flow \eqref{eq:reduced-slow-sub} is singular along the fold $\mathcal{L_S}$ \eqref{eq:fold}. We desingularize via the transformation $t_s = \det(J)t_d$, which gives the desingularized reduced system on $\mathcal{M_S}$
\begin{subequations}\label{eq:desin-reduced-slow-sub}
\begin{align}
    \frac{dx_1}{dt_d} &= (x_2^2-1) g_1(x_1,x_2)=F_1\\
   \frac{dx_2}{dt_d} &= (x_1^2-1) g_2(x_1,x_2) = F_2
\end{align}
\end{subequations}
The desingularized system is equivalent to \eqref{eq:reduced-slow-sub} when $\det(J)>0$ but has the opposite orientation otherwise. 
The folded singularity $M$ of the desingularized reduced system \eqref{eq:desin-reduced-slow-sub} are isolated points that lie on one of the four fold curves and at which the right hand side of \eqref{eq:desin-reduced-slow-sub} is zero. These are the special points that can allow trajectories to cross the fold with nonzero speed. For example, in the stable folded node case, trajectories that land inside some trapping region on $\mathcal{M_S}$ will converge to the folded node, thereby passing through the fold from an attracting branch to a repelling branch. Such solutions are so-called singular canards. We refer the readers to \cite{Desroches2012} for review. 

\begin{remark}\label{rm:fs}
A folded singularity in $M$ is 
either $\{x_1=\pm 1\}$ (i.e., $\mathcal{L}_{1\pm}$) and $F_1=0$ or $\{x_2=\pm 1\}$ (i.e., $\mathcal{L}_{2\pm}$) and $F_2=0$.
\end{remark}

\begin{figure}[htp!]
    \centering   
    \includegraphics[width=3.5in]{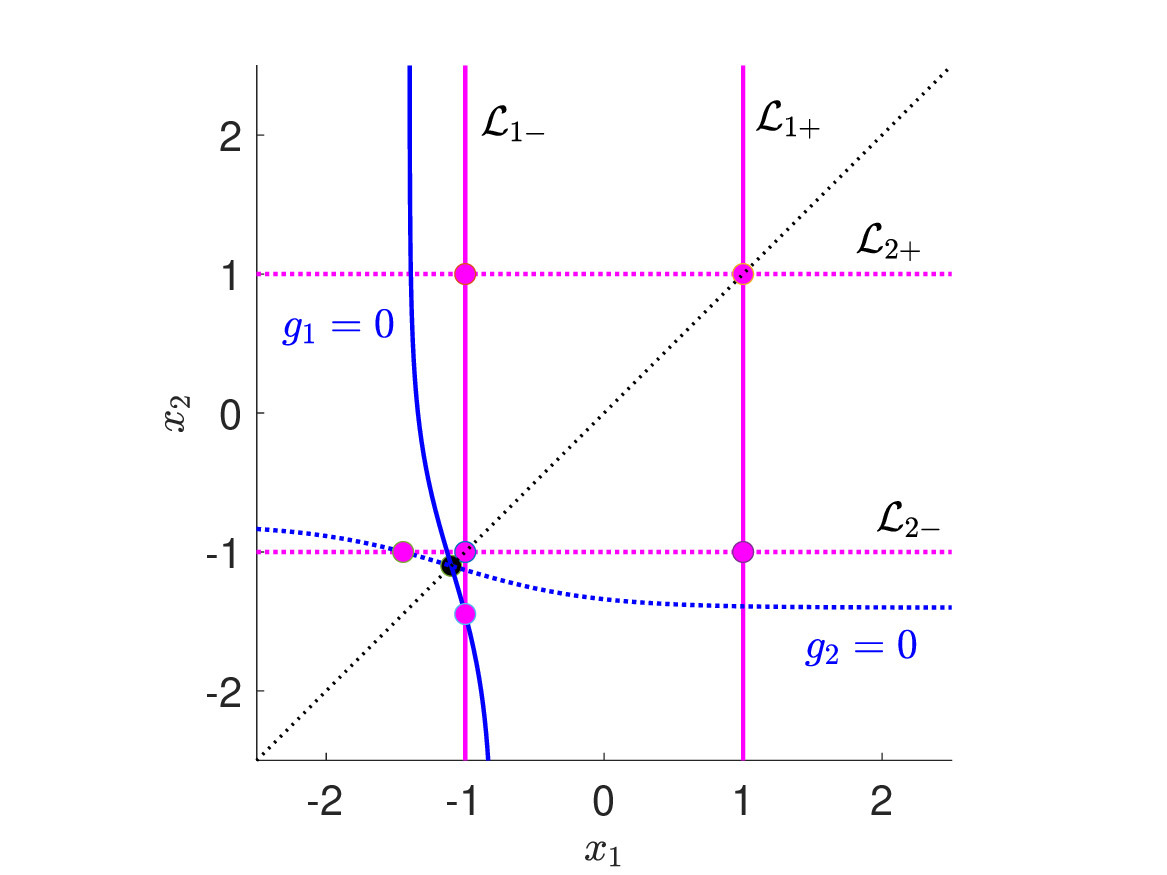}
    \caption{Fold curves (magenta lines), $y_i$-nullclines (blue lines) and folded singularity points (magenta circles) of the coupled van der Pol oscillators \eqref{eq:cVdP}, projected onto $(x_1,x_2)$-space. Parameter values are given in Table \ref{tab:parameters}. The black circle denotes the synchronous full system equilibrium $(-1.1, -1.1)$. Folded singularities near the true equilibrium are folded saddles. Those further away are stable folded focus (($(-1, 1)$ and $(1, -1)$)) or folded saddle ($(1,1)$). 
    }
    \label{fig:fs}
\end{figure}

Figure \ref{fig:fs} displays the four fold curves (magenta solid curves for $x_1$-nullcline folds; dashed for $x_2$-nullcline folds), $g_1=0$ ($y_1$-nullcline, blue solid), and $g_2=0$ ($y_2$ nullcline, blue dashed). Note that a folded singularity is not a true equilibrium of the full system \eqref{eq:cVdP}, i.e., the intersection of $y_1$- and $y_2$- nullclines denoted by the black circle in Figure \ref{fig:fs}. It follows that there are a total of six folded singularities, denoted by magenta circles in Figure \ref{fig:fs}. 
Specifically, the folded singularities of \eqref{eq:cVdP} are given by 
\begin{itemize}
    \item Along $\mathcal{L}_{1-}$ ($x_1=-1$), $x_2^2-1=0 \rightarrow x_2=\pm 1$ or $g_1=0 \rightarrow x_2 = -1.4466$.
    \item Along $\mathcal{L}_{1+}$ ($x_1=1$), $x_2=\pm 1$ (note that $g_1\neq0$ along $\mathcal{L}_{1+}$)
    \item Along $\mathcal{L}_{2-}$ ($x_2=-1$), $x_1^2-1=0 \rightarrow x_1=\pm 1$ or $g_2=0 \rightarrow x_1 = -1.4466$.
    \item Along $\mathcal{L}_{2+}$, $x_2=1$, $x_2=\pm 1$ (note that $g_2\neq0$ along $\mathcal{L}_{2+}$)
\end{itemize} 

The eigenvalues of the Jacobian of the desingularized system evaluated along $M$ are used to classify the folded singularites. Folded singularities with two real eigenvalues with the same sign (resp., with opposite signs) are called \textit{folded nodes} (resp., \textit{folded saddles}). Those with complex eigenvalues are called \textit{folded foci}, which does not produce canard dynamics. In our system with parameters given in Table \ref{tab:parameters}, the three folded singularity points near the true equilibrium (the black circle in Figure \ref{fig:fs}) are all folded saddles, so is the folded singularity $(1,1)$. The remaining two folded singularities, $(-1, 1)$ and $(1, -1)$, are folded foci. 

In the coupled oscillators studied in \cite{awal2024symmetry}, the two folded singularity points near the true equilibrium, located off the symmetry axis, are both folded nodes. These folded nodes are the key mechanism driving the strong symmetry breaking in their system. Specifically, passage through the neighborhood of a folded node results in splitting between the amplitudes of the two oscillations. This mechanism of symmetry-breaking is not present in our system. 

\bibliographystyle{unsrt}
\bibliography{references}

\end{document}